\documentclass[11pt]{article}
\usepackage{amsfonts,mathrsfs,amssymb,amsthm,mathptm}
\usepackage{amsmath,amscd}
\usepackage{mathptm,pslatex}
\usepackage{color}
\usepackage[dvips]{graphicx}
\usepackage[all]{xy}
\usepackage{graphicx}
\usepackage{fancyhdr}

\begin{document}
\newtheorem{Def}{Definition}[section]
\newtheorem{Bsp}[Def]{Example}
\newtheorem{Prop}[Def]{Proposition}
\newtheorem{Theo}[Def]{Theorem}
\newtheorem{Lem}[Def]{Lemma}
\newtheorem{Koro}[Def]{Corollary}
\theoremstyle{definition}
\newtheorem{Rem}[Def]{Remark}

\newcommand{\add}{{\rm add}}
\newcommand{\con}{{\rm con}}
\newcommand{\gd}{{\rm gl.dim}}
\newcommand{\sd}{{\rm st.dim}}
\newcommand{\sr}{{\rm sr}}
\newcommand{\dm}{{\rm dom.dim}}
\newcommand{\cdm}{{\rm codomdim}}
\newcommand{\tdim}{{\rm dim}}
\newcommand{\E}{{\rm E}}
\newcommand{\Mor}{{\rm Morph}}
\newcommand{\End}{{\rm End}}
\newcommand{\ind}{{\rm ind}}
\newcommand{\rsd}{{\rm res.dim}}
\newcommand{\rd} {{\rm rd}}
\newcommand{\ol}{\overline}
\newcommand{\overpr}{$\hfill\square$}
\newcommand{\rad}{{\rm rad}}
\newcommand{\soc}{{\rm soc}}
\renewcommand{\top}{{\rm top}}
\newcommand{\pd}{{\rm pdim}}
\newcommand{\id}{{\rm idim}}
\newcommand{\fld}{{\rm fdim}}
\newcommand{\Fac}{{\rm Fac}}
\newcommand{\Gen}{{\rm Gen}}
\newcommand{\fd} {{\rm fin.dim}}
\newcommand{\Fd} {{\rm Fin.dim}}
\newcommand{\Pf}[1]{{\mathscr P}^{<\infty}(#1)}
\newcommand{\DTr}{{\rm DTr}}
\newcommand{\cpx}[1]{#1^{\bullet}}
\newcommand{\D}[1]{{\mathscr D}(#1)}
\newcommand{\Dz}[1]{{\mathscr D}^+(#1)}
\newcommand{\Df}[1]{{\mathscr D}^-(#1)}
\newcommand{\Db}[1]{{\mathscr D}^b(#1)}
\newcommand{\C}[1]{{\mathscr C}(#1)}
\newcommand{\Cz}[1]{{\mathscr C}^+(#1)}
\newcommand{\Cf}[1]{{\mathscr C}^-(#1)}
\newcommand{\Cb}[1]{{\mathscr C}^b(#1)}
\newcommand{\Dc}[1]{{\mathscr D}^c(#1)}
\newcommand{\K}[1]{{\mathscr K}(#1)}
\newcommand{\Kz}[1]{{\mathscr K}^+(#1)}
\newcommand{\Kf}[1]{{\mathscr  K}^-(#1)}
\newcommand{\Kb}[1]{{\mathscr K}^b(#1)}
\newcommand{\DF}[1]{{\mathscr D}_F(#1)}

\newcommand{\Kac}[1]{{\mathscr K}_{\rm ac}(#1)}
\newcommand{\Keac}[1]{{\mathscr K}_{\mbox{\rm e-ac}}(#1)}

\newcommand{\modcat}{\ensuremath{\mbox{{\rm -mod}}}}
\newcommand{\Modcat}{\ensuremath{\mbox{{\rm -Mod}}}}
\newcommand{\pres}{\ensuremath{\mbox{{\rm pres}}}}
\newcommand{\Spec}{{\rm Spec}}

\newcommand{\stmc}[1]{#1\mbox{{\rm -{\underline{mod}}}}}
\newcommand{\Stmc}[1]{#1\mbox{{\rm -{\underline{Mod}}}}}
\newcommand{\prj}[1]{#1\mbox{{\rm -proj}}}
\newcommand{\inj}[1]{#1\mbox{{\rm -inj}}}
\newcommand{\Prj}[1]{#1\mbox{{\rm -Proj}}}
\newcommand{\Inj}[1]{#1\mbox{{\rm -Inj}}}
\newcommand{\PI}[1]{#1\mbox{{\rm -Prinj}}}
\newcommand{\GP}[1]{#1\mbox{{\rm -GProj}}}
\newcommand{\GI}[1]{#1\mbox{{\rm -GInj}}}
\newcommand{\gp}[1]{#1\mbox{{\rm -Gproj}}}
\newcommand{\gi}[1]{#1\mbox{{\rm -Ginj}}}

\newcommand{\opp}{^{\rm op}}
\newcommand{\otimesL}{\otimes^{\rm\mathbb L}}
\newcommand{\rHom}{{\rm\mathbb R}{\rm Hom}\,}
\newcommand{\pdim}{\pd}
\newcommand{\Hom}{{\rm Hom}}
\newcommand{\Coker}{{\rm Coker}}
\newcommand{ \Ker  }{{\rm Ker}}
\newcommand{ \Cone }{{\rm Con}}
\newcommand{ \Img  }{{\rm Im}}
\newcommand{\Ext}{{\rm Ext}}
\newcommand{\StHom}{{\rm \underline{Hom}}}
\newcommand{\StEnd}{{\rm \underline{End}}}
\newcommand{\KK}{I\!\!K}
\newcommand{\gm}{{\rm _{\Gamma_M}}}
\newcommand{\gmr}{{\rm _{\Gamma_M^R}}}

\def\vez{\varepsilon}\def\bz{\bigoplus}  \def\sz {\oplus}
\def\epa{\xrightarrow} \def\inja{\hookrightarrow}

\newcommand{\lra}{\longrightarrow}
\newcommand{\llra}{\longleftarrow}
\newcommand{\lraf}[1]{\stackrel{#1}{\lra}}
\newcommand{\llaf}[1]{\stackrel{#1}{\llra}}
\newcommand{\ra}{\rightarrow}
\newcommand{\dk}{{\rm dim_{_{k}}}}

\newcommand{\holim}{{\rm Holim}}
\newcommand{\hocolim}{{\rm Hocolim}}
\newcommand{\colim}{{\rm colim\, }}
\newcommand{\limt}{{\rm lim\, }}
\newcommand{\Add}{{\rm Add }}
\newcommand{\Prod}{{\rm Prod }}
\newcommand{\Tor}{{\rm Tor}}
\newcommand{\Cogen}{{\rm Cogen}}
\newcommand{\Tria}{{\rm Tria}}
\newcommand{\Loc}{{\rm Loc}}
\newcommand{\Coloc}{{\rm Coloc}}
\newcommand{\tria}{{\rm tria}}
\newcommand{\Con}{{\rm Con}}
\newcommand{\Thick}{{\rm Thick}}
\newcommand{\thick}{{\rm thick}}
\newcommand{\Sum}{{\rm Sum}}

{\Large \bf
\begin{center}
Derived equivalences, new matrix equivalences, and homological conjectures
\end{center}}

\medskip
\centerline{\textbf{Xiaogang Li} and \textbf{Changchang Xi}$^*$ }

\renewcommand{\thefootnote}{\alph{footnote}}
\setcounter{footnote}{-1} \footnote{ $^*$ Corresponding author.
Email: xicc@cnu.edu.cn; Fax: 0086 10 68903637.}
\renewcommand{\thefootnote}{\alph{footnote}}
\setcounter{footnote}{-1}
\renewcommand{\thefootnote}{\alph{footnote}}
\setcounter{footnote}{-1}

\begin{abstract} Based on the fact that every finite-dimensional algebra over a field is isomorphic to the centralizer of \textbf{two} matrices, we approach the representation theory of finite-dimensional algebras over fields by centralizers of matrices. The first fundamental question is to study the centralizer of a single matrix, called a centralizer matrix algebra. By introducing three new equivalence relations on all square matrices over a field, we completely characterize Morita, derived and almost $\nu$-stable derived equivalences between centralizer matrix algebras in terms of these matrix equivalences, respectively. Further, we show that a derived equivalence between centralizer matrix algebras of permutation matrices induces both a Morita equivalence and additional derived equivalences for $p$-regular parts and for $p$-singular parts. As an application, we show that the finitistic dimension conjecture and the Nakayama conjecture are valid for centralizer matrix algebras.
\end{abstract}

{\footnotesize\tableofcontents\label{contents}}

\section{Introduction\label{Introduction}}

Derived categories and equivalences between them are the pi\`ece de r\'esistance of modern
homological algebra. They were initiated by Grothendieck around 1960's and developed further by Verdier (see \cite{verdier}). Since then a lot of applications and connections have been discovered to other branches in mathematics. For instance, in representation theory, Happel applied them successfully to generalized tilting modules over finite-dimensional algebras \cite{happel88}. Moreover, Rickard advanced Happel's work and developed a beautiful Morita theory
for derived categories of rings (see \cite{Rickard1, JR2}). Also, Keller established a Morita theory for differential graded
algebras (see \cite{keller}). All of these provide powerful tools to understand derived module categories and
equivalences of both rings and differential graded rings. However, it is still a hard and untractable, but fundamental, problem to decide whether two algebras are derived equivalent or not. This can be seen from a not yet solved conjecture by Brou\'e, which says that a block algebra of a finite group algebra with abelian defect subgroup should be derived equivalent to its Brauer corresponding block algebra \cite{broue1990}. Though many efforts have been made in the last decades, the conjecture seems far away from being solved completely. For some advances about this conjecture, we refer to \cite{chuang-rouquier, rouq}.

Based on the fact that every finite-dimensional algebra over a field is isomorphic to the centralizer algebra of \textbf{two} matrices \cite[Lemma 2]{brenner1972}, we try to understand general finite-dimensional noncommutative algebras by centralizer algebras of matrices.  This point of view is completely different from the present situation in the representation theory of algebras by quivers with relations. We start first with a fundamental step by studying the centralizer algebra of a single matrix, called a \emph{centralizer matrix algebra}.
The purpose of this article is:

(1) To introduce new equivalence relations for all square matrices and to provide complete descriptions of Morita, derived and almost $\nu$-stable derived equivalences for  centralizer matrix algebras.

(2) To show that the Nakayama conjecture and the finitistic dimension conjecture are valid for centralizer matrix algebras over fields.

\smallskip
Let us state our main results and their consequences more precisely.

Let $R$ be a field. We denote by $M_n(R)$ the full $n\times n$ matrix algebra over $R$ with the identity matrix $I_n$.
For a nonempty subset $X$ of $M_n(R)$, the centralizer algebra $S_n(X,R)$ of $X$ in $M_n(R)$ is defined by $$S_n(X,R):=\{a\in M_n(R)\mid ax=xa,\; \forall \; x\in X\}.$$

Clearly, $S_n(X,R)=\cap_{c\in X}S_n(\{c\},R)$. For simplicity,  we write $S_n(c,R)$ for $S_n(\{c\},R)$, and term $S_n(c,R)$ as a \emph{centralizer matrix algebra} in this article.

It is easy to see that any finite-dimensional algebra $A$ over a field $R$ is isomorphic to the centralizer algebra of finitely many matrices in $M_n(R)$ with $n=\dim_R(A)$. But, for a nonempty finite set $X$, $S_n(X,R)$ is isomorphic to $S_m(Y,R)$ for some $m\in \mathbb{N}$ and $Y$ with $|Y|= 2$ by \cite[Lemma 1]{brenner1972}. Thus every finite-dimensional algebra over a field can be realized as the centralizer algebra of two matrices (see also \cite[Lamma 2]{brenner1972}). This implies that the study of centralizer matrix algebras is an important ingredient of the representation theory and homological algebra of finite-dimensional algebras.

On centralizer matrix algebras, Frobenius proved a nice dimension formula in terms of the degrees of  invariant factors of the given matrix (see \cite[Theorem 1, Theorem 2, p.105-106]{We1934}).

\medskip
\textbf{Theorem} (Frobenius). Let $d_1(x),\cdots, d_s(x)$ be the invariant factors of positive degree of a matrix $c\in M_n(R)$ over a field $R$, and let $n_i$ be the degree of $d_i(x), 1\le i\le s$. Then  $\dim_RS_n(c,R)=\sum_{i=1}^s(2s-2i+1)n_i$. \label{frob-dim}

Centralizer matrix algebras include centrosymmetric matrix algebras (see \cite{W85, xy}) and the quasi-hereditary  Auslander algebras of the truncated polynomial algebras $R[x]/(x^n)$ for all $n$ (see \cite {xz1}), which play a crucial role in the classification of  parabolic subgroups of classical groups with a finite number of orbits on the unipotent radical (see \cite{hr}). Also, all algebras of the form $R[x]/(f(x))$ can be realized as centralizer matrix algebras.
If $c$ is an invertible matrix, then the centralizer matrix algebra of $c$ is the invariant algebra of the  action of cyclic group $\langle c\rangle$ on $M_n(R)$ by conjugation.
In general, if $X$ consists of invertible matrices, then $S_n(X,R)$ is just the invariant algebra which dates back to the classical invariant theory (see \cite{HW}). If $c$ is a nilpotent matrix in $M_n(\mathbb{F}_q)$, where $\mathbb{F}_q$ is a finite field with $q$ elements, then the determinants of matrices in $S_n(c,\mathbb{F}_q)$ are completely described (see \cite{bw}). Centralizer matrix algebras are also studied in invariant orbits (see \cite{bf}), and in maximal doubly stochastic matrix theory (see \cite{cdfk}).

Recently, a lot of new structural and homological properties of $S_n(c,R)$ have been revealed in a series of papers \cite{xz1,xz2,xz3}. For instance, $S_n(c,R)$ is always a cellular $R$-algebra in the sense of Graham--Lehrer (see \cite{gl1996}) if the field $R$ is algebraically closed. Moreover, the famous Auslander--Reiten conjecture (or Auslander--Alperin conjecture) on stable equivalences holds true for centralizer matrix algebras \cite{xz3}. Recall that the conjecture states that stably equivalent algebras should have the same number of non-isomorphic, non-projective simple modules \cite[Conjecture (5), p.409]{ARS} (see also \cite{rouq}).

Concerning centralizer matrix algebras, there are still many basic problems. In this article, we focus on the following question.

\smallskip
{\bf Question}: Let $R$ be a field, $c\in M_n(R)$ and $d\in M_m(R)$.  What are necessary and sufficient conditions for $S_n(c,R)$ and $S_m(d,R)$ to be Morita or derived equivalent?

\smallskip
To answer this question, we introduce the so-called $M$-equivalence, $D$-equivalence and $AD$-equivalence on all square matrices over fields. These matrix equivalences reflect information on maximal elementary divisors of matrices. We refer the reader to Section \ref{sect3} for precise definitions.
A complete answer to the above question reads as follows.

\begin{Theo}\label{main1}
Let $R$ be a field, $c\in M_n(R)$ and $d\in  M_m(R)$. Then
the centralizer matrix algebras of $c$ and of $d$ are Morita equivalent (respectively, derived equivalent, or almost $\nu$-stable derived equivalent) if and only if the matrices $c$ and $d$ are $M$-equivalent (respectively, $D$-equivalent, or $AD$-equivalent).
\end{Theo}

Thus the existence of a Morita equivalence, a derived equivalence or an almost $\nu$-stable derived equivalence between centralizer matrix algebras can be read off directly from the elementary divisors of given matrices, and therefore is reduced to matrix equivalences in linear algebra.

As an application of our methods, we consider  the Nakayama conjecture \cite{Nakayama} and the finitistic dimension conjecture \cite{bass}.

\medskip
\textbf{Nakayama Conjecture} (NC): An Artin algebra is self-injective if it has infinite dominant dimension.

\smallskip
\textbf{Finitistic Dimension Conjecture} (FDC): The finitistic dimension of an Artin algebra is always finite.

\smallskip
These are two of the central conjectures in the representation theory and homological algebra of Artin algebras (see \cite[Conjectures, p.409]{ARS}). They are still open up to date. But we will show in Subsection \ref{sect4.2} that the conjectures hold true for centralizer matrix algebras.

\begin{Theo}\label{fdc}
$(1)$ The finitistic dimension conjecture holds true for centralizer matrix algebras over fields. Particularly, the Nakayama conjecture holds true for centralizer matrix algebras over fields.

$(2)$ If two centralizer matrix algebras are derived equivalent, then they have the same dominant dimension.
\end{Theo}

Consequently, (FDC) is valid for any algebras that are derived equivalent to centralizer matrix algebras because the finiteness of finitistic dimensions is invariant of derived equivalences.

Next, we state some corollaries of Theorem \ref{main1}. For unexplained notation, we refer to Subsection \ref{sect3.1}.

\begin{Koro}\label{derp}Let $R$ be a field, $c\in M_n(R)$ and $d\in  M_m(R)$.

$(1)$ If $c$ and $d$ are permutation matrices, then $S_n(c,R)$ and $S_m(d,R)$ are Morita equivalent if and only if they are derived equivalent.

$(2)$ If the field $R$ is perfect, then the following  are equivalent.

$\quad (a)$ $S_n(c,R)$ and $S_m(d,R)$ are almost $\nu$-stable derived equivalent.

$\quad (b)$ $S_n(c,R)$ and $S_m(d,R)$ are stably equivalent of Morita type, and there is a bijection $\pi: \mathcal{M}_c\setminus\mathcal{R}_c\ra \mathcal{M}_d\setminus\mathcal{R}_d$, such that $R[x]/(f(x))\simeq R[x]/((f(x))\pi)$ for $f(x)\in \mathcal{M}_c\setminus\mathcal{R}_c$.
\end{Koro}

For a derived equivalence of the centralizer matrix algebras of permutation matrices, we can additionally get two more derived equivalences for $p$-regular and $p$-singular parts of the permutations, where $p$ is the characteristic of the ground field. The $p$-regular part $r(\sigma)$ and the $p$-singular part $s(\sigma)$ of $\sigma \in \Sigma_n$ are defined in terms of the cycle type of $\sigma$. For more details, we refer to Subsection \ref{sect4}.

Let $c_{\sigma}:=\sum_{i=1}^ne_{i,(i)\sigma}\in M_n(R)$ be the permutation matrix of $\sigma$, where $e_{ij}$ is the matrix with $1$ in $(i,j)$-entry and $0$ in all other entries.

\begin{Koro}[Proposition \ref{regular-singular}]\label{cor1.5}
Let $R$ be a field of characteristic $p\ge 0$, $\sigma\in \Sigma_n$ and $\tau\in \Sigma_m$. If $S_n(c_{\sigma},R)$ and $S_m(c_{\tau},R)$ are derived equivalent, then

$(1)$ $S_n(c_{r(\sigma)},R)$ and $S_m(c_{r(\tau)},R)$ are derived equivalent, and

 $(2)$ $S_n(c_{s(\sigma)},R)$ and $S_m(c_{s(\tau)},R)$ are derived equivalent.
\end{Koro}

The paper is organized as follows. In Section \ref{sect2} we fix notation, recall basic definitions and terminology, and prove a few preliminary lemmas needed in the later proofs. In Section \ref{sect3} we introduce 3 new equivalence relations on square matrices over fields. As examples, we describe representation-finite centralizer matrix algebras. In Section \ref{Pf} we prove the main results and their corollaries. In Section \ref{examples} we present examples to show that the converse of Corollary \ref{cor1.5} may be false and that even for centralizer matrix algebras over a field, the notions of  Morita, derived and almost $\nu$-stable derived equivalences are distinct, though they may coincide in many cases. Finally, we propose some open problems for further investigation. For example, can one generalize the main results for $R$ being a principal ideal domain?

\section{Preliminaries}\label{sect2}

In this section we recall some basic definitions and terminologies on derived equivalences, and prepare a few lemmas on modules over polynomial algebras for our proofs.

\subsection{Definitions and notation\label{sect2.1}}
In this paper, $R$ is a field unless stated otherwise. By an algebra we mean a finite-dimensional unitary associative algebra over $R$. By a module we mean a left module.

Let $A$ be an algebra. By $\rad(A)$ and $LL(A)$ we denote the Jacobson radical and Loewy length of $A$, respectively. Let $A^{\opp}$ and $A^e$ stand for the opposite algebra and the enveloping algebra $A\otimes_R A\opp$ of $A$, respectively.

We write $A\modcat$ for the category of all finitely generated $A$-modules, $A\modcat_{\mathscr{P}}$ for the full subcategory of $A\modcat$ consisting of modules without any nonzero projective direct summands, and $A\prj$ (respectively, $A\inj$) for the full subcategory of $A\modcat$ consisting of projective (respectively, injective) $A$-modules.

For an $A$-module $M\in A\modcat$, $\ell(M)$ denotes the composition length of $M$, and ${\add}(M)$ denotes the full subcategory of $A\modcat$ consisting of all modules isomorphic to direct summands of direct sums of finitely many copies of $M$. If $M\in A\prj$, we denote by $\pres(M)$ the full subcategory of $A\modcat$ consisting of those modules $L$ such that there is an exact sequence $P_1\ra P_0\ra L\ra 0$ with $P_0, P_1\in \add(M)$. The\emph{ basic module} of $M$ is by definition the direct sum of all non-isomorphic indecomposable direct summands of $M$. This is uniquely determined by $M$ up to isomorphism, and denoted by $\mathcal{B}(M)$. Let $M_{\mathscr{P}}$ be the submodule of $M$ such that $M_{\mathscr{P}}$ has no nonzero projective direct summand and $M/M\!_{\mathscr{P}}$ is projective. Thus $M\!_{\mathscr{P}}\in A\modcat\!_{\mathscr{P}}$.

For homomorphisms $f:X\to Y$ and $g: Y\to Z$ in $A\modcat$, we write $fg$ for their composition. This implies that the image of an element $x\in X$ under $f$ is denoted by $(x)f$. Thus $\Hom_A(X,Y)$ is naturally an $\End_A(X)$-$\End_A(Y)$-bimodule, where $\End_A(X)$ stands for the endomorphism algebra of the module $X$.

The composition of functors between categories is written from right to left, that is, for two functors
$F:\mathcal{C}\ra \mathcal{D}$ and $G:\mathcal{D}\ra \Sigma$, we write $G\circ F$, or simply $GF$, for the composition of $F$ with $G$. The image of an object $X\in \mathcal{C}$ under $F$ is written as $F(X).$

Let $\mathcal{D}$ be a class of $A$-modules. By the number of modules in $\mathcal{D}$ we always mean the number of the isomorphism classes of modules in $\mathcal{D}$.

A homomorphism $f: M\ra N$ in $A\modcat$ is  \emph{right almost split} if $f$ is not a split surjection and any homomorphism $X\ra N$ which is not a split surjection factorizes through $f$. Dually, left almost split homomorphisms are defined. An exact sequence $0\ra M\stackrel{f}\ra L\stackrel{g}\ra N\ra 0$ of $A$-modules is called an \emph{almost split sequence} if $f$ is left almost split and $g$ is right almost split. We refer to \cite{ARS} for further information on almost split sequences.
The homomorphism $f$ is called a \emph{radical homomorphism} if, for any $Z\in A\modcat$, $g\in \Hom_A(Z,M)$ and $h\in\Hom_A(N,Z)$, the composition $gfh$ is not an automorphism of $Z$.

Let $D=\Hom_R(-,R): A\modcat\to A^{\opp}\modcat$ be the usual duality of $A$. The Nakayama functor $\nu_A:=D\Hom_A(-,A)\simeq D(A)\otimes_A-: A\modcat \to A\modcat$ restricts to an equivalence between $\prj{A}$ and $A\inj$. An $A$-module $M\in A\modcat$ is said to be $\nu$-\emph{stably projective} if $\nu^i_A{M}$ is projective for all $i\ge 0$. Let $A$-stp denote the full subcategory of $A\modcat$ consisting of all $\nu$-stably projective $A$-modules.
Clearly, there is an idempotent $e\in A$ such that $A$-stp = $\add(Ae)$. The self-injective algebra $eAe$ is called the \emph{Frobenius part} of $A$, which is unique up to Morita equivalence (see \cite{hx2} or \cite{MV2} for more details).

The $R$-algebra $A$ is said to be \emph{elementary} if $A/\rad(A)$ is isomorphic to the direct product of copies of $R$, and  \emph{split} if there exist positive integers $n_1,\cdots,n_s$ such that $A/\rad(A)\simeq \oplus^{s}_{j=1}M_{n_j}(R)$ as algebras. So elementary $R$-algebras are always split.

Let $\Db{A}$ stand for the bounded derived category of $A\modcat$. It is known that $\Db{A}$ is an $R$-linear, triangulated category.
Let $\stmc{A}$ denote the stable module category of $A\modcat$,
which is the quotient category of $A\modcat$ modulo the full subcategory $A$-proj. In general, $\stmc{A}$ is not a triangulated category. But, if $A$ is self-injective, then $\stmc{A}$ is an $R$-linear triangulated category.

\begin{Def}
Algebras $A$ and $B$ over a field $R$ are said to be

$(1)$ Morita equivalent if their module categories $A\modcat$ and $B\modcat$ are equivalent as $R$-linear categories. In this case, an equivalence $F: A\modcat\rightarrow B\modcat$ of $R$-linear categories is called a Morita equivalence between $A$ and $B$.

$(2)$ Derived equivalent if their derived categories $\Db{A}$ and $\Db{B}$ are equivalent as $R$-linear triangulated categories. In this case, an equivalence $F:\Db{A}\rightarrow\Db{B}$ of $R$-linear triangle categories is called a derived equivalence between $A$ and $B$.

$(3)$ Stably equivalent if their stable module categories $A\stmc$ and $B\stmc$ are equivalent as $R$-linear categories. In this case, an equivalence $F: \stmc{A} \ra  B\stmc$ of $R$-linear categories is called a stable equivalence between $A$ and $B$.
\end{Def}

For further information on derived categories and equivalences of rings, we refer to \cite{Rickard1,JR2}.

If $F$ is a stable equivalence between algebras $A$ and $B$, then $F$ induces a one-to-one correspondence between non-isomorphic, indecomposable, non-projective modules in $A\modcat_{\mathscr{P}}$ and  $B\modcat_{\mathscr{P}}$.

The following is a simple observation on Morita equivalences.

\begin{Lem}\label{add}
Let $A$ be an algebra and $M,N\in A\modcat$. Then ${\End}_A(M)$ and $\End_A(N)$ are Morita equivalent if and only if $\add(M)$ and $\add(N)$ are equivalent as $R$-linear categories. Moreover, if $A$ is a local, Nakayama algebra, then the algebras $\End_A(M)$ and $\End_A(N)$ are Morita equivalent if and only if the basic modules $\mathcal{B}(M)$ and $\mathcal{B}(N)$ are isomorphic.
\end{Lem}

{\it Proof.} We only prove the second statement.
Suppose that $A$ is a local, Nakayama algebra with $LL(A)=n$. Then $\{A/\rad^i(A)\mid 0\le i \le n-1\}$ is a complete list of all non-isomorphic indecomposable $A$-modules, and $\End_A(A/\rad^i(A))\simeq A/\rad^i(A)$ as algebras for $0\leq i\leq n-1$. Thus,  for indecomposable $A$-modules $X$ and $Y$, $\End_A(X)\simeq \End_A(Y)$ if and only if $X\simeq Y$. Suppose that ${\End}_A(M)$ and $\End_A(N)$ are Morita equivalent. Then there is an $R$-linear equivalence $G: {\add}(M)\ra {\add}(N)$. In particular, we have $\End_A(C)\simeq \End_A(G(C))$ for  $C\in {\add}(M)$, and therefore $C\simeq G(C)$ for any indecomposable module $C\in \add(M)$. Hence $\mathcal{B}(M)\simeq \mathcal{B}(N)$ as $A$-modules. The converse is clear by the fact:  $\add(M)=\add(\mathcal{B}(M))$ for any $M\in A\modcat$. $\square$

The second statement in Lemma \ref{add} is not true in general. For example, if $A$ is an algebra over an algebraically closed field $R$ and has at least two (non-isomorphic) simple $A$-modules $M$ and $N$ such that $\End_A(M)\simeq R\simeq \End_A(N)$, then we cannot get $M\simeq N$.

\medskip
As a special class of derived equivalences, almost $\nu$-stable derived equivalences were introduced in \cite{hx1}. Recall that a tilting complex is called a \emph{radical tilting complex} if all of its differentials are radical homomorphisms. Every tilting complex over an algebra $A$ is isomorphic to a radical tilting complex in  $\Db A$ (see \cite[(a), p.112]{hx1}).

\begin{Def}{\rm \cite{hx1}} Let $F:\Db{A}\ra\Db{B}$ be a derived equivalence of algebras $A$ and $B$. Suppose that $\cpx{Q}$ and $\cpx{\bar{Q}}$ are radical tilting complexes associated to $F$ and the quasi-inverse $F^{-1}$ of $F$, respectively. By applying the shift functor if necessary, we may assume that $\cpx{Q}$ and $\cpx{\bar{Q}}$ are of the form
$$0\lra Q^{-n}\lra\cdots \lra Q^{-1}\lra Q^0\lra 0,\quad 0\lra {\bar Q}^0\lra {\bar Q}^1\lra\cdots\lra {\bar Q}^n\lra 0,$$
respectively. Let $Q:=\bigoplus_{i=1}^nQ^{-i}$ and $\bar{Q}:=\bigoplus_{i=1}^n\bar{Q}^n$. The derived equivalence $F$ is said to be {\em almost $\nu$-stable} provided that $\add({}_AQ)=\add(\nu_AQ)$ and $\add({}_B\bar{Q})=\add(\nu_B\bar{Q})$.
\end{Def}

One of the significant properties of  almost $\nu$-stable derived equivalences is that such an equivalence between algebras always induces a stable equivalence of Morita type (see \cite[Theorem 1.1]{hx1}), and thus preserves global and dominant dimensions of algebras. This generalises a result of Rickard on derived equivalences of self-injective algebras (see \cite[Corollary 5.5]{JR2}).

\begin{Def} {\rm \cite{MB}} Algebras $A$ and $B$ are \emph{stably equivalent of Morita type} if there exist bimodules $_AM_B$ and $_BN_A$ such that $M$ and $N$ are projective as one-sided modules, $M\otimes_B N\simeq A\oplus P$ and $N\otimes_A M\simeq B\oplus Q$ as bimodules, where $P$ is a projective $A^e$-module and $Q$ is a projective $B^e$-module.
\end{Def}

In this definition, the exact functor $N\otimes_A-: A\modcat{}\ra B\modcat$ induces a stable equivalence $N\otimes_A-:A\stmc{}\ra B\stmc$.

Examples of stable equivalences of Morita type are the derived equivalences between self-injective algebras (see \cite[Corollary 5.5]{JR2}). Another example is that a commutative ring $R$ and a separable $R$-algebra $A$ are stably equivalent of Morita type. Here, an $R$-algebra $A$ is \emph{separable} over $R$ if $_AA_A$ is a projective $A^e$-module.

\smallskip
An algebra is said to be \emph{representation-finite} if it has only finitely many non-isomorphic indecomposable modules. Consequently, given $A$-modules $M$ and $N$ with $\add(N)\subseteq \add(M)$, if $\End_A(M)$ is representation-finite, then so is $\End_A(N)$. Equivalently, if $A$ is representation-finite, then so is $eAe$ for all $e=e^2\in A$.

\subsection{Basic facts on derived equivalences of algebras}
Derived equivalences of algebras were described by Rickard in terms of tilting complexes in \cite{Rickard1}. However, for our purpose, we will follow the approach in \cite{hx2} to construct derived equivalences of algebras. For further information on constructing derived equivalences of algebras, we refer to \cite{x3}.

Let $\mathcal{C}$ be an additive category and $\mathcal{D}$ a full subcategory of $\mathcal{C}$. Given an object $Y\in \mathcal{C},$ a morphism $f:M\ra Y$ in $\mathcal{C}$ is called a \emph{right $\mathcal{D}$-approximation} of $Y$ if $M\in \mathcal{D}$ and each morphism $D\to Y$ with $D\in \mathcal{D}$ factorizes through $f$. A \emph{left $\mathcal{D}$-approximation} of an object $X$ in $\mathcal{C}$ is defined dually. As usual, we denote by $\End_{\mathcal{C}}(Y)$ the endomorphism ring of an object $Y\in\mathcal{C}$.

\begin{Def} {\rm \cite{hx2}} A sequence $X\stackrel{g}\ra M\stackrel{f}\ra Y$ of morphisms in $\mathcal{C}$ with $M\in \mathcal{D}$ is called a $\mathcal{D}$-split sequence if $g$ is both a kernel of $f$ and a left $\mathcal{D}$-approximation of $X$, and if $f$ is both a cokernel of $g$ and a right $\mathcal{D}$-approximation of $Y$.
\end{Def}

Examples of ${\add}(M)$-split sequences capture almost split sequences $X\to M\to Z$ in $A\modcat$. Also, for any projective-injective module $M$ and a submodule $X$ of $M$, the exact sequence $X\ra M\ra M/X$ is an $\add(M)$-split sequence.

\begin{Lem}{\rm\cite[Theorem~ 1.1]{hx2}}
\label{split-thm}
Let $A$ be an algebra, and let $\mathcal{C}$ be a full additive subcategory of $A\modcat$ and $M$  an object in $\mathcal{C}$. Suppose
that $X\ra M'\ra Y$ is an ${\add}(M)$-split sequence in $\mathcal{C}$. Then ${\End}_{\mathcal{C}}(M\oplus X)$ and
 ${\End}_{\mathcal{C}}(M\oplus Y)$ are derived equivalent.
\end{Lem}

As a consequence of Lemma \ref{split-thm}, we get the following result (see also \cite[Section 3, Remark]{hx1}).
\begin{Lem}\label{alm}
Let $A$ be a self-injective algebra and $X\in A\modcat$. Then ${\End}_A(A\oplus X)$ and ${\End}_A(A\oplus \Omega_A(X))$ are almost $\nu$-stable derived equivalent.
\end{Lem}

The next result is somehow a converse of Lemma \ref{alm}.
\begin{Lem}\label{almst}
{\rm \cite[Theorem 4.4]{CM}}
Let $A$ and $B$ be symmetric algebras, and let $F$ be an almost $\nu$-stable derived equivalence between ${\rm End}_A(A\oplus M)$ and ${\rm End}_B(B\oplus N)$, where $_AM$ and $_BN$ are basic non-zero modules without nonzero projective summands. Then $A$ and $B$ are (almost $\nu$-stable) derived equivalent. Furthermore, $F$ induces a stable equivalence $\overline{F}:A\stmc$ $\ra B\stmc$ with $\overline{F}(M) = N$.
\end{Lem}

\begin{Lem}\label{iso}
Let $A$ and $B$ be commutative self-injective algebras, and let $_AM$ and $_BN$ be faithful modules over $A$ and $B$, respectively. If ${\End}_A(M)$ and ${\End}_B(N)$ are derived equivalent, then $A\simeq Z(\End_A(M))\simeq Z\big({\End}_B(N))\simeq B$, where $Z(C)$ denotes the center of an algebra $C$.
\end{Lem}

{\it Proof.} For an algebra $C$ and a faithful $C$-module $X$, one always has an embedding $Z(C)\hookrightarrow Z(\End_C(X))$. Thus $A\hookrightarrow Z(\End_A(M))$ since $A$ is commutative. Note that a faithful module over a self-injective algebra is clearly a generator-cogenerator. This implies that $M_{\End_A(M)}$ is a right faithful module and the bimodule $_AM_{\End_A(M)}$ has the double centralizer property, that is $\End_{\End_A(M)\opp}(M)\simeq A$. Thus there is an embedding $Z(\End_A(M))\hookrightarrow \End_{\End_A(M)\opp}(M)\simeq A.$ Hence $A\simeq Z(\End_A(M))$.
Now, assume that ${\End}_A(M)$ and ${\End}_B(N)$ are derived equivalent. Then $Z({\End}_A(M)\simeq Z({\End}_B(N))$ by \cite[Proposition 9.2]{Rickard1}, and therefore $A\simeq Z(\End_A(M))\simeq Z({\End}_B(N))\simeq B$. $\square$

\subsection{Modules over quotients of polynomial algebras\label{2.3}}
In this subsection we recall some basic facts on modules over the polynomial algebra $R[x]$, where $R$ is a field, and prove a few basic lemmas for later proofs.

Throughout this section, let $f(x)$ be a fixed irreducible polynomial in $R[x]$ and $A:=R[x]/(f(x)^n)$ for a natural number $n>0$. Then $A$ is a local, commutative, symmetric, Nakayama algebra (see, for instance \cite[Example, p.127]{ARS}). Thus $A$ has $n$ indecomposable modules $M(i):=R[x]/(f(x)^i)$ for $i\in [n].$ We write $M(0)=0$. Clearly, $\Hom_A(M(i),A)\simeq \Hom_R(M(i),R)\simeq M(i)$ as $A$-modules, and $\ell(M(i))=i$ for all $i\in [n]$. Moreover,
for $i,j\in [n]$, we see that $i\le j$ if and only if there is an injective homomorphism in $\Hom_A(M(i),M(j))$ if and only if there is a surjective homomorphism in $\Hom_A(M(j),M(i))$.

For $B:=R[x]/(f(x)^m)$ with $m<n$, there is a canonical surjective homomorphism $\pi:A\ra B$ of $R$-algebras, and therefore each $B$-module can be viewed as an $A$-module via $\pi$. Up to isomorphism, indecomposable $A$-modules coming from $B$-modules are exactly those $M(i)$ with $i\in [m].$ Clearly, ${\Hom}_A(M,N)={\Hom}_B(M, N)$ for $M, N\in B\modcat$.

For an irreducible polynomial $g(x)\in R[x]$ and a positive integer $m$, if $A\simeq R[x]/(g(x)^m)$ as $R$-algebras, then $n=LL\big(R[x]/(f(x)^n)\big)=LL\big(R[x]/(g(x)^m)\big)=m$ and, for $t\in [n]$, the indecomposable $B$-module $R[x]/(g(x)^t)$ is isomorphic to the $A$-module $M(t)$.

Now, suppose that $G:\stmc{A}  \to A\stmc$ is a stable equivalence. For $n\ge 2$, we define $\Gamma_{n-1}:=\{M(i)\mid i\in [n-1]\}\subseteq A\modcat_{\mathscr{P}}$. Then $G$ induces a permutation $\ol{G}$ on $\Gamma_{n-1}$, namely, for  $M\in \Gamma_{n-1}$, $\ol{G}(M)$ is the unique module in $\Gamma_{n-1}$ such that $\ol{G}(M)\simeq G(M)$ in $A\stmc$. Clearly, $\overline{\Omega}_A(M(i))= M(n-i)$, where $\Omega_A$ is the syzygy operator of $A$.

\begin{Lem}\label{self} Let $n\ge 2$. If $G:\stmc{A}  \to A\stmc$ is a stable equivalence, then the induced action $\ol{G}$ on $\Gamma_{n-1}$ coincides with either $\overline{\Omega}_A$ or the identity action.
\end{Lem}

{\it Proof.} If $n=2$, then $A$ has only one non-projective indecomposable $A$-module $S$ and $\Omega(S)\simeq S$. Thus the conclusion is true. Now suppose $n\geq 3$.  For $A$-modules $X$ and $Y$, let Irr$(X,Y)$ be the $R$-space $\rad_A(X,Y)/\rad_A^2(X,Y)$. For $i,j\in [n-1]$, it follows from the shape of the Auslander-Reiten quivers of Nakayama algebras that Irr$(M(i),M(j))\neq 0$ if and only if $|i-j|\leq 1$. By a general result on stable equivalences (see \cite[Lemma 1.2, p.336]{ARS}), we have Irr$(X,Y)\simeq$Irr$(G(X),G(Y))$ as $R$-spaces for $X,Y\in A\modcat_{\mathscr{P}}.$ It then follows that $G(M(1))\simeq M(1)$ or $G(M(1))\simeq M(n-1)=\Omega_A(M(1))$.  This implies that $G(M(i))\simeq M(i)$ for $i\in [n-1]$ or $G(M(i))\simeq \Omega_A(M(i))$ for $i\in [n-1]$. Hence $\ol{G}$ is the identity map or equals $\overline{\Omega}_A$. $\square$

\begin{Lem}\label{lift-exact}
Let $a,b,c,d \in \{0,1,\cdots,n\}$ such that $b < a< c$, $b< d< c$ and $a+d=b+c$. If $_AX\in A\modcat$ has no indecomposable direct summands $N$ with $b<\ell(N)<c$ and $_AY:={}_AX\oplus M(b)\oplus M(c)$,  then there is an ${\add}(_AY)$-split sequence $0\ra M(a)\ra M(b)\oplus M(c)\ra M(d)\ra 0$.
\end{Lem}

{\it Proof.}
Let $g:M(b)\to M(d)$ and $h: M(c)\to M(d)$ be the canonical injective and surjective homomorphisms, respectively, and define
$v:=\left(\begin{smallmatrix} g\\	h\\		\end{smallmatrix} \right).$
Then $v:M (b)\oplus M(c)\to M(d)$ is a surjective homomorphism. Similarly, let $p: M(a)\to M(b)$ and $q: M(a)\to M(c)$ be the canonical surjective and injective homomorphisms, respectively, and define $u:=(-p,q)$. Then $u: M(a)\to M(b)\oplus M(c)$ is an injective homomorphism. By the definition of $M(i)$, we have $uv=0.$ It follows from $a+d=b+c$ that the sequence
$$ (*)\quad 0\lra M(a)\stackrel{u}\lra M(b)\oplus M(c)\lraf{v} M(d)\lra 0$$
of $A$-modules is exact. This can also be seen  from the Auslander--Reiten quivers of Nakayama algebras.

We shall show that $u$ and $v$ are left and right ${\add}(_AY)$-approximations of $M(a)$ and $M(d)$, respectively. In fact, we need only to show that $v$ is a right ${\add}(_AY)$-approximation of $M(d)$ because the dual functor ${\Hom}_R(-,R)$ transforms right ${\add}(_AY)$-approximations to left ${\add}(_AY)$-approximations. To show that $v$ is a right ${\add}(_AY)$-approximation of $M(d)$, it suffices to prove that, for any indecomposable direct summand $Z$ of $_AY$, each homomorphism $h:Z\to M(d)$ factorizes through $v$. Let $\Img(h)$ denote the image of $h$. By the assumption on $X$, we have either $\ell(Z)\leq b$ or $\ell(Z)\geq c.$

Suppose $\ell(Z)\leq b.$ Then $\ell(\Img(h))\leq \ell(Z)\leq b=\ell((M(b))g)$. Since $(M(b))g$ is maximal submodule of $M(d)$ of length $b$,  we have $\Img(h)\subseteq (M(b))g.$ Let $s: Z\to M(b)\oplus M(c)$ be the map defined by $(z)s:=((z)h)g^{-1},0)$ for $z\in Z.$ Clearly, $h$ is a homomorphism of $A$-modules such that $h=sv$, that is, $h$ factorizes through $v$.

Suppose $\ell(Z)\geq c.$ Then max$\{a,b,c, d\}\le\ell(Z)$. Let $B:=R[x]/(f(x)^{\ell(Z)})$. Then $B$ is the quotient of $A$ by the ideal $(f(x))^{n-\ell(Z)}+(f(x)^n)$ and $Z\simeq M(\ell(Z))=R[x]/(f(x)^{\ell(Z)})=B$ as $A$-modules. This shows that $Z$ is also a projective $B$-module. Thus the exact sequence $(*)$ can be viewed as a sequence of $B$-modules. So the exactness of ${\Hom}_{B}(Z,-)$ implies that $h$ factorizes through $v$ in $B\modcat$. Since ${\Hom}_A(M, N)={\Hom}_B(M,N)$ for $M, N\in B\modcat$, we see that $h$ factorizes through $v$ in $A\modcat$. $\square$

\begin{Lem}\label{de}
Let $n=\sum^s_{i=1} \ell_i$ with $\ell_i\in \mathbb{Z}_{>0}$. For $\sigma\in \Sigma_s$ and $j\in [s]$, define $M_j:=M(\sum^j_{i=1} \ell _i)$ and  $M^\sigma_j:= M(\sum^j_{i=1} \ell_{(i)\sigma})$. Then $\End_A(\bigoplus^s_{j=1} M_j)$ and ${\End}_A(\bigoplus^s_{j=1} M^\sigma_j)$ are derived equivalent.
\end{Lem}

{\it Proof.} The symmetric group $\Sigma_s$ is generated by the transpositions $(t,t+1),t\in [s-1]$. In particular, $\sigma\in \Sigma_s$ can be written as a product of these transpositions, say $\sigma=\prod^k_{i=1} (t_i,t_i+1)$ for $t_i\in [s-1].$ Set $\sigma_{k+1}:=id$ and $\sigma_r:=\prod^k_{i=r} (t_i,t_i+1)$ for all $r\in [k].$ Then $\sigma_1=\sigma$ and $(t_r,t_r+1)\sigma_r=\sigma_{r+1}$ for $r\in [k]$. In particular, $(t_r)\sigma_{r+1}=(t_r+1)\sigma_r, (t_r+1)\sigma_{r+1}=(t_r)\sigma_r$ and $(t)\sigma_{r+1}=(t)\sigma_r$ for $t\in [s]\setminus \{t_r,t_r+1\}$.

Since $ \End_A(\bigoplus^s_{j=1} M^\sigma_j) = {\End}_A\big(\bigoplus^s_{j=1} M^{\sigma_1}_j\big)$ and $\End_A(\bigoplus^s_{j=1} M_j)=\End_A\big(\bigoplus^s_{j=1} M^{\sigma_{k+1}}_j  $, it suffices to show that there is a derived equivalence between $\End_A\big(\bigoplus^s_{j=1} M^{\sigma_r}_j\big)$ and $\End_A\big(\bigoplus^s_{j=1} M^{\sigma_{r+1}}_j\big)$ for all $r\in [k].$

Indeed, for any $\tau\in \Sigma_s$, we define $\sum^{t_r-1}_{i=1} \ell_{(i)\tau}=0$ if $t_r=1$. For $r\in [k]$, let $a_r:=\ell_{(t_r+1)\sigma_{r+1}}+\sum^{t_r-1}_{i=1} \ell_{(i)\sigma_{r+1}},$ $b_r:=\sum^{t_r-1}_{i=1} \ell_{(i)\sigma_{r+1}},c_r:=\sum^{t_r+1}_{i=1} \ell_{(i)\sigma_{r+1}},d_r:=\sum^{t_r}_{i=1} \ell_{(i)\sigma_{r+1}}$, $X_r:=\bigoplus_{j\in [s],|j-t_r|\geq 2} M^{\sigma_{r+1}}_j $ and $Y_r:=\bigoplus_{j\in [s], j\neq t_r} M^{\sigma_{r+1}}_j$. Then $b_r<a_r<c_r,b_r<d_r<c_r, a_r+d_r=b_r+c_r$ and
$$Y_r = M\big(\sum^{t_r-1}_{i=1} \ell_{(i)\sigma_{r+1}}\big)\oplus M\big(\sum^{t_r+1}_{i=1} \ell_{(i)\sigma_{r+1}}\big) \oplus \bigoplus_{j\in [s],|j-t_r|\geq 2} M^{\sigma_{r+1}}_j= M(b_r)\oplus M(c_r)\oplus X_r.$$
Clearly, for any indecomposable direct summand $Z$ of $X_r=\bigoplus_{j\in [s],|j-t_r|\geq 2} M^{\sigma_{r+1}}_j$, either $\ell(Z)\leq \sum^{t_r-2}_{j=1}\ell_{(j)\sigma_{r+1}}$ $ < b_r$ or $\ell(Z)\geq \sum^{t_r+2}_{j=1}\ell_{(j)\sigma_{r+1}}>c_r$.
It then follows from Lemma \ref{lift-exact} that there is an ${\add}(Y_r)$-split sequence
$$0\lra M\big(a_r)\lra M(b_r)\oplus M(c_r)\lra M\big(d_r)\lra 0.$$ Clearly, $\bigoplus^s_{j=1} M^{\sigma_r}_j= Y_r\oplus M(a_r)$ and $\bigoplus^s_{j=1} M^{\sigma_{r+1}}_j=Y_r\oplus M(d_r).$ By Lemma \ref{split-thm}, ${\End}_A(\bigoplus^s_{j=1} M^{\sigma_r}_j)$ and ${\End}_A(\bigoplus^s_{j=1} M^{\sigma_{r+1}}_j)$ are derived equivalent. $\square$

\begin{Rem}\label{mult-eq}
The sums $\sum^j_{i=1} \ell _i$ and $\sum^j_{i=1} \ell_{(i)\sigma}$, appearing in Lemma \ref{de}, are related to the definition of $D$-equivalences of matrices (see Section \ref{newequrel} below). For $s\geq 2$ and a series of integers $m_s>m_{s-1}>\cdots>m_1\geq 1$, let $\ell_1:=m_1$ and $\ell_i:=m_i-m_{i-1}$ for $2\leq i\leq s$. Then $m_j=\sum^j_{i=1}\ell_i$ for $j\in [s].$ For another series of integers $n_s>n_{s-1}>\cdots>n_1\geq 1$, if  $\{\{m_s-m_{s-1},\cdots,m_1\}\}=\{\{n_s-n_{s-1},\cdots,n_1\}\}$ as multisets, then there exists some $\sigma\in \Sigma_s$ such that $n_j=\sum^j_{i=1}\ell_{(i)\sigma}$ for $j\in [s].$ Moreover, if $\{\{m_s-m_{s-1},\cdots,m_1\}\}=\{\{n_s-n_{s-1},\cdots,n_1\}\}$ and if there are two irreducible polynomials $f(x)$ and $g(x)$ in $R[x]$ such that $R[x]/(f(x)^{m_s})\simeq R[x]/(g(x)^{n_s})$ as algebras, then it follows from Lemma \ref{de} that $\End_{R[x]/(f(x)^{m_s})}\big(\bigoplus_{k\in [s]} R[x]/(f(x)^{m_k})\big)$ and $\End_{R[x]/(g(x)^{n_s})}\big(\bigoplus_{k\in [s]} R[x]/(g(x)^{n_k})\big)$ are derived equivalent.
\end{Rem}

Recall that a polynomial $g(x)\in R[x]$ of positive degree is \emph{separable} if it has only simple roots in its splitting field.
\begin{Lem} \label{poi}
If the irreducible polynomial $f(x)$ is separable, then $K:=R[x]/(f(x))$ is a separable field over $R$, the algebra $A$ can be viewed as a $K$-algebra, and $A\simeq K[x]/(x^n)$ as $K$-algebras.
\end{Lem}

{\it Proof.} Since $f(x)$ is separable and $\rad(A)=(f(x))/(f(x)^n)$, we know that $A/\rad(A)\simeq K$ is a separable $R$-algebra. By Wedderburn-Malcev Theorem \cite[Theorems 24 and 28]{We}, there exists a subalgebra $S$ of $A$ such that $A=S\oplus \rad(A)$ as $R$-vector spaces. Consequently, $S\simeq A/\rad(A)\simeq K.$ So $A$ can be viewed as a $K$-algebra. Since $A$ is a finite-dimensional, elementary, local $K$-algebra of representation-finite type, there is an $m\in \mathbb{N}$ such that $A\simeq K[x]/(x^m)$ as $K$-algebras. By considering the chain $R[x]\supsetneq (f(x))\supsetneq (f(x)^2)\supsetneq\cdots\supsetneq (f(x)^n)\supsetneq 0$ and comparing the $K$-dimensions of the algebras in the isomorphism, we get $n = m$. $\square$

\begin{Koro} \label{St-i}
If the polynomial $f(x)$ is separable and $g(x)\in R[x]$ is irreducible such that $A$ is stably equivalent to $R[x]/(g(x)^m)$ for an integer $m\ge 2$, then $A\simeq R[x]/(g(x)^m)$ as $R$-algebras and $m=n$.
\end{Koro}

{\it Proof.}  Since stably equivalent algebras of representation-finite type have the same number of non-isomorphic, non-projective, indecomposable modules, we have $n-1=m-1$, and therefore $n=m$. Set $B:=R[x]/(g(x)^m)$. Let $F: A\stmc{}\ra B\stmc$ be a stable equivalence, and let $S$ be the unique simple $A$-module (up to isomorphism). Then $F$ induces a one-to-one correspondence between the set of non-isomorphic, non-projective, indecomposable modules in $A\modcat_{\mathscr{P}}$ and the one in
$B\modcat_{\mathscr{P}}$. Thanks to $n=m\ge 2$, the module $S$ is not projective and $\End_A(S)\simeq \underline{\End}_A(S)$. Thus $F(S)$ is indecomposable and $$\End_A(S)\simeq \underline{\End}_A(S)\simeq \underline{\End}_B(F(S))=\End_B(F(S))/\mathcal{P}(F(S),F(S))$$is a division ring, where $\mathcal{P}(F(S),F(S))$ is the set of all homomorphisms that factorize through projective $B$-modules.
Since $\mathcal{P}(F(S),F(S))\subseteq \rad(\End_B(F(S)))$, we have $\mathcal{P}(F(S),F(S))=\rad(\End_B(F(S)))$. This yields the following isomorphisms of algebras:
$$R[x]/(f(x))\simeq A/\rad(A)\simeq\underline{\End}_{A}(S)\simeq \underline{\End}_B(F(S))\simeq B/\rad(B)\simeq R[x]/(g(x)).$$
In particular, $g(x)$ is also a separable polynomial. Let $K:=R[x]/(f(x))$. Then Lemma \ref{poi} implies that $A\simeq K[x]/(x^n)\simeq B$ as $K$-algebras, and therefore also as $R$-algebras. $\square$

\medskip
For $c\in M_n(R)$,
set $A_c:=R[x]/(\Ker(\varphi))=R[x]/(m_c(x)) \simeq R[c]$. Then the characteristic matrix $xI_n-c$ of $c$ is a matrix over the principal ideal domain $R[x]$. Suppose that  $xI_n-c$ has invariant factors $d_1(x),\cdots, d_r(x)$ of positive degree with $r\le n$ and $d_i|d_{i+1}$ for $1\le i<r$. Let $d_r(x)=f_1(x)^{e_{r1}}\cdots f_s(x)^{e_{r s}}$, where $f_1(x),\cdots, f_s(x)$ are pairwise coprime, irreducible polynomials, and $e_{rj}>0$ is an integer for $j\in [s]$. Then, for $i\in [r-1]$, we can write $d_i(x)=f_1(x)^{e_{i1}}\cdots f_s(x)^{e_{i s}}$, where $0\le e_{i j}\le e_{i+1 j}\le \cdots\le e_{rj}$. The polynomials $f_j(x)^{e_{ij}}$, with $e_{ij}$ positive for $ i\in [r]$ and $j\in [s]$, are called the \emph{elementary divisors} of $c$. This can be interpreted alternatively in the following way.

Let $R^n$ be the set of $n\times 1$ matrices with entries in $R$. Then $c$ can be viewed as a linear transformation $\sigma_c$ on $R^n$ by ${\sigma_c}\cdot v:=cv$ for $v\in R^n$. Note that $m_c(x)=d_r(x)=f_1(x)^{e_{r1}}\cdots f_s(x)^{e_{r s}}$. Set $M_j:=\Ker(f_j(\sigma_c)^{e_{rj}})$ for $j\in [s]$. Then $M_j$ is a $\sigma_c$-invariant subspace (equivalently, $R[c]$-submodule) of $R^n$ and
$ R^n= \bigoplus^s_{j=1}M_j$.
Note that the minimal polynomial of the restriction of $\sigma_c$ to $M_j$ is $f_j(x)^{e_{rj}}$. By \cite[Theorem 4.11]{Co}, we see that $M_j=\oplus^{l_j}_{i=1}M_{ji}$ can decompose into direct sum of $\sigma_c$-cyclic subspaces $M_{ji}$, and the minimal polynomial of the restriction of $\sigma_c$ to the subspace $M_{ji}$ is $f_j(x)^{q_{ji}}$. The multiset of these polynomials $f_j(x)^{q_{ji}}$ is in fact the multiset of elementary divisors of $c$ (over $R$). Note that $R^n$ can be regarded as an $R[x]$-module by letting $x^k$ act on $R^n$ as $\sigma^k_c$ and the decomposition $R^n= \bigoplus^s_{j=1}\bigoplus^{l_j}_{i=1} \; M_{ji}$ is in fact a decomposition of $R^n$ as a direct sum of indecomposable submodules. That the minimal polynomial of the restriction of $\sigma_c$ to the subspace $M_{ji}$ is $f_j(x)^{q_{ji}}$ is equivalent to saying that $M_{ji}\simeq R[x]/(f_j(x)^{q_{ji}})$ as $R[x]$-modules. Thus
$$(\star)\quad \quad R^n\simeq \bigoplus^s_{j=1}\bigoplus^{l_j}_{i=1} \; R[x]/(f_j(x)^{q_{ji}})$$
as $R[x]$-modules (see \cite[Chapter 4, p.130-133]{Co} for more details).

For $c\in M_n(R)$, we deine $\mathcal{E}_c$ to be the set of elementary divisors of $c$. Here, we understand that a set always has no duplicate elements. Further, we define the set of \emph{maximal divisors} of $c$ by
$$\mathcal{M}_c:=\{f(x)\in \mathcal{E}_c\mid f(x) \mbox{ is maximal with respect to polynomial divisibility} \}.$$

The next lemma follows immediately from $(\star)$.

\begin{Lem} \label{bijection}  If $R^n$ is identified with the $A_c$-module $\bigoplus^s_{j=1}\bigoplus^{l_j}_{i=1}R[x]/(f_j(x)^{q_{ji}})$ in $(\star)$, then there is a bijection $\pi$ from $\mathcal{E}_c$ to the set of pairwise non-isomorphic indecomposable direct summands of the $A_c$-module $R^n$, sending $h(x)$ to the $A_c$-module $R[x]/(h(x))$ for $h(x)\in \mathcal{E}_c$.
\end{Lem}

Suppose that the characteristic of $R$ is  $p\geq 0$. For a positive integer $m$, there exist uniquely determined integers $s, m'\in \mathbb{N}$ such that $m=p^s m'$ and $p\nmid m'$, we define $\nu_p(m):=s$. Here, we understand $\nu_p(m):=0$ if $p=0$.
Suppose that $\sigma\in \Sigma_n$ is a permutation of cycle type $(\lambda_1,\cdots,\lambda_k)$. Let $g(x)$ be an irreducible factor of the minimal polynomial $m_{c_\sigma}(x)$ of the permutation matrix $c_{\sigma}$ of $\sigma$, we define $q_{g(x)}:= max\{\nu_{p}(\lambda_j) \mid j\in[k], \mbox{ such that } g(x)\mbox{ divides } x^{\lambda_j}-1  \}$.  Note that $q_{g(x)}$ depends upon the cycle type of $\sigma$.

\begin{Lem}\label{per}
Suppose that the characteristic of $R$ is $p\ge 0$ and $\sigma\in \Sigma_n$ is a permutation of cycle type $(\lambda_1,\cdots,\lambda_k)$. Then $$\mathcal{E}_{c_{\sigma}} = \{g(x)^{p^{\nu_p(\lambda_i)}}\mid i\in[k], g(x)~\mbox{is an irreducible factor of}~ x^{\lambda_i}-1 \} \mbox{ and  }$$
$$\mathcal{M}_{c_{\sigma}} = \{g(x)^{p^{q_{g(x)}}}\mid  g(x)~\mbox{is an irreducible factor of}~ m_{c_{\sigma}}(x) \}.$$
\end{Lem}

{\it Proof.} For conjugate permutations in $\Sigma_n$, their corresponding permutation matrices are similar, and therefore have the same elementary divisors. Thus, without loss of generality, we may assume that $\sigma=(1,\cdots,\lambda_1)(\lambda_1+1,\cdots,\lambda_1+\lambda_2)\cdots(\sum^{k-1}_{j=1}\lambda_j+1,\cdots,n)$. Then $c_\sigma$ is a diagonal block matrix, that is $c_\sigma$= diag$\{c_{\sigma_1}, c_{\sigma_2}, \cdots, c_{\sigma_k}\}$, where $\sigma_i$ is a $\lambda_i$-cycle in $\Sigma_{\lambda_i}$ for $i\in [k]$.  In particular, $m_{c_\sigma}(x)$ is the least common multiple of $x^{\lambda_i}-1$ for $i\in [k]$ and $\mathcal{E}_{c_{\sigma}} =\bigcup_{i\in [k]}\mathcal{E}_{c_{\sigma_i}}$. For a matrix $d\in M_m(R)$, let $\chi_d(x)$ denote the characteristic polynomial of $d$ in $R[x]$. For $i\in [k]$, we write $\lambda_i=p^{\nu_p(\lambda_i)}\lambda_i'$ with $p\nmid \lambda'_i$. Then $x^{\lambda'_i}-1=\prod^{h_i}_{j=1}f_{ij}(x)$, where  $f_{i1}(x), f_{i2}(x), \cdots, f_{ih_i}(x)$ are distinct irreducible (monic) polynomials in $R[x]$. From the following equalities
$$\chi_{c_{\sigma_i}}(x)=x^{\lambda_i}-1=x^{p^{\nu_p(\lambda_i)}\lambda_i'}-1 =(x^{\lambda'_i}-1)^{p^{\nu_p(\lambda_i)}}=
\prod^{h_i}_{j=1}f_{ij}(x)^{p^{\nu_p(\lambda_i)}},$$ we get
$\chi_{c_{\sigma_i}}(x)=m_{c_{\sigma_i}}(x)=x^{\lambda_i}-1$. Hence $\mathcal{E}_{c_{\sigma_i}}=\mathcal{M}_{c_{\sigma_i}}$. This implies $$\mathcal{E}_{c_{\sigma}} = \{g(x)^{p^{\nu_p(\lambda_i)}}\mid i\in[k], g(x)~\mbox{is an irreducible factor of}~ x^{\lambda_i}-1\}.$$

Clearly, $\mathcal{M}_{c_\sigma}$ is of the form $\{g_1(x)^{m_1},g_2(x)^{m_2},\cdots, g_t(x)^{m_t}\}$, where $g_1(x), g_2(x),$ $\cdots, g_t(x)$ form a complete set of distinct (monic) irreducible factors of $m_{c_\sigma}(x)$ and where $m_1,m_2,\cdots, m_t$ are positive integers. Let  $g_s(x)$ be an irreducible factor of $m_{c_\sigma}(x)$. Then $g_s(x)$ divides $x^{\lambda_i}-1$ for at least one $i\in [k]$, and therefore the set $S(g_s(x)):=\{g_s(x)^{p^{\nu_p(\lambda_j)}}\mid j\in [k] \mbox{ and }~g_s(x)~\mbox{divides}~x^{\lambda_j}-1\}\neq \varnothing$. By the description of $\mathcal{E}_{c_\sigma}$, we see that $S(g_s(x))$ is exactly the elementary divisors of $c_\sigma$ which are divided by $g_s(x)$. Thus $g_s(x)^{p^{q_{g_s(x)}}}$ is a maximal elementary divisor of $c_{\sigma}$ by the definition of $q_{g_s(x)}$. Hence $$\mathcal{M}_{c_{\sigma}} = \{g(x)^{p^{q_{g(x)}}}\mid  g(x)~\mbox{is an irreducible factor of}~ m_{c_{\sigma}}(x) \}. \quad \square$$

\smallskip
Now, we prove a result on congruences of matrices that appear as the Cartan matrices of the endomorphism algebras of modules over polynomial algebras.
Note that two multisets $\{\{x_1,\cdots,x_s\}\}$ and $\{\{y_1,\cdots,y_s\}\}$ are equal if and only if there exists a permutation $\sigma\in \Sigma_s$ such that $(y_1,\cdots,y_s)^\sigma:=(y_{(1)\sigma},\cdots,y_{(s)\sigma})=(x_1,\cdots,x_s).$

\begin{Lem}\label{mat}
For an integer $s\geq 2$, let $m_1> m_2> \cdots >m_s\ge 1$ and $n_1>n_2>\cdots>n_s\ge 1$ be two series of integers with $m_1=n_1$. Set $X:=\sum^s_{k=1} \big(\sum^k_{l=1} m_k(e_{kl}+e_{lk})-m_k e_{kk}\big)\in M_s(\mathbb{Z})$ and $Y:=\sum^s_{k=1} \big(\sum^k_{l=1} n_k(e_{kl}+e_{lk})-n_k e_{kk}\big)\in M_s(\mathbb{Z})$. Then $X$ and $Y$ are congruent in $M_s(\mathbb{Z})$ if and only if there is $\sigma\in \Sigma_s$ such that $(n_1-n_2,\cdots,n_{s-1}-n_s,n_s)=(m_1-m_2,\cdots,m_{s-1}-m_s,m_s)^\sigma$.
\end{Lem}

{\it Proof.}  We define three matrices in $M_s(\mathbb{Z})$ by $U:=I_s-\sum^{s-1}_{t=1}e_{t,t+1},$ $D_1:={\rm diag}(m_1-m_2,\cdots,m_{s-1}-m_s,m_s) \mbox{ and } D_2 :={\rm diag}(n_1-n_2,\cdots,n_{s-1}-n_s,n_s).$ Then $U^{tr}XU = D_1$ and $U^{tr} Y U = D_2$, where $U^{tr}$ stands for the transpose of $U$. Thus $X$ and
$Y$ are congruent in $M_s(\mathbb{Z})$ if and only if $D_1$ and $D_2$ are congruent in $M_s(\mathbb{Z}).$ Now, we show that $D_1$ and $D_2$ are congruent in $M_s(\mathbb{Z})$ if and only if there is an element $\sigma\in \Sigma_s$ such that $(n_1-n_2,\cdots,n_{s-1}-n_s,n_s)=(m_1-m_2,\cdots,m_{s-1}-m_s,m_s)^\sigma.$ Indeed, if $(n_1-n_2,\cdots,n_{s-1}-n_s,n_s)=(m_1-m_2,\cdots,m_{s-1}-m_s,m_s)^\sigma$ for some $\sigma\in \Sigma_s$, then $c^{tr}_{\sigma} D_1 c_{\sigma}=D_2$. This means that $D_1$ and $D_2$ are congruent in $M_s(\mathbb{Z}).$ Conversely, suppose that $D_1$ and $D_2$ are congruent in $M_s(\mathbb{Z})$. Then there is an invertible matrix $H=(h_{ij})_{1\leq i,j\leq s}\in M_s(\mathbb{Z})$ such that $H^{tr} D_1 H=D_2.$ This implies\vspace{-0.3cm}
$$(*)\quad \sum^{s-1}_{r=1} (\sum^s_{k=1}h^2_{kr})(m_r-m_{r+1})+(\sum^s_{k=1}h^2_{ks})m_s=n_1=m_1.\vspace{-0.2cm}$$
Since $H$ is invertible in $M_s(\mathbb{Z})$, each column of $H$ has a nonzero element, and therefore $\sum^s_{k=1}h^2_{kr}\ge 1$ for $r\in [s].$ Now it follows from ($*$) that  $\sum^s_{k=1}h^2_{kr}=1$ for all $r\in [s]$. Thus each row and column of $H$ has only one nonzero entry which is either $1$ or $-1$. This implies that $H = \epsilon c_{\tau}$ for $\tau\in \Sigma_s$ and a diagonal matrix $\epsilon$ with the entries in $\{1,-1\}$.  Hence $H^{tr}=H^{-1}.$ This shows that the diagonal matrices $D_1$ and $D_2$ are similar, and therefore they have the same eigenvalues (counting multiplicities). So $\{\{m_1-m_2,\cdots,m_{s-1}-m_s,m_s\}\}= \{\{n_1-n_2,\cdots,n_{s-1}-n_s,n_s\}\}$ as multisets, that is, $(n_1-n_2,\cdots,n_{s-1}-n_s,n_s)=(m_1-m_2,\cdots,m_{s-1}-m_s,m_s)^\sigma$ for some $\sigma\in \Sigma_s.$ $\square$

\section{New equivalence relations of matrices\label{sect3}}
In this section we introduce three new equivalence relations on square matrices over a field, and present necessary and sufficient conditions for centralizer matrix algebras to be representation-finite.

\subsection{Definitions of matrix equivalences\label{sect3.1}}

Let $R[x]$ be the polynomial algebra over a field $R$ in one variable $x$.
Given polynomials $f(x)$ and  $g(x)$ of positive degree, if $f(x)$ divides $g(x)$, that is, $g(x)=f(x)h(x)$ with $h(x)\in R[x]$, we write $f(x)\mid g(x)$.  Observe that this divisibility of polynomials defines a partial order on the set of all monic polynomials of positive degree in $R[x]$.

Let $n$ be a natural number and $c\in M_n(R)$. Recall that $\mathcal{E}_c$ denotes the set of elementary divisors of $c$, and
$\mathcal{M}_c:=\{f(x)\in \mathcal{E}_c\mid f(x) \mbox{ is maximal with respect to polynomial divisibility} \}$ is called the set of maximal divisors of $c$. In fact, $\mathcal{M}_c$ is determined completely by the invariant factor $d_r(x)$ or $m_c(x)$.

Let $\mathcal{R}_c:=\{f(x)\in \mathcal{M}_c\mid f(x) \mbox{ is reducible}\}$. This is the set of all reducible maximal divisors of $c$.

For $f(x)\in \mathcal{M}_c$, we define the set $P_c(f(x))$ of \emph{power indices} in $\mathcal{E}_c$ by

$P_c(f(x)):=\{i\ge 1 \mid  \exists \mbox{ irreducible polynomial } p(x) \mbox{ such that } p(x) \mbox{ divides } f(x), p(x)^i\in \mathcal{E}_c \}.$

Let $\mathbb{Z}_{>0}$ be the set of all positive integers and $s\in \mathbb{Z}_{>0}$. For a subset  $T:=\{ m_1, m_2, \cdots, m_s\}$ of $\mathbb{Z}_{>0}$ with $m_1>m_2>\cdots >m_s$, we define a set $\mathcal{J}_T:=\{m_1,m_1-m_2,\cdots,m_1-m_s\}$ and a \emph{multiset} $\mathcal{H}_T:= \{\{m_1-m_2,\cdots,m_{s-1}-m_s,m_s\}\}$. Note that we allow duplicate elements to occur in multisets. If $s=1$, then $\mathcal{H}_T=\mathcal{J}_T=T$. Observe that if $H=\{n_1, n_2, \cdots, n_s\}$ is another subset of $\mathbb{Z}_{>0}$ with $n_1>n_2>\cdots>n_s$, then $H=\mathcal{J}_T$ if and only if $T=\mathcal{J}_H$.

\smallskip
Now we introduce three new equivalence relations on the set of all square matrices over a field.

\begin{Def}\label{newequrel}
Two matrices $c\in M_n(R)$ and $d\in M_m(R)$ are said to be

$(1)$ \emph{$M$-equivalent} if there is a bijection $\pi:\mathcal{M}_c\ra\mathcal{M}_d$, such that $R[x]/(f(x)) \simeq R[x]/((f(x))\pi)$ as algebras and ${P_c(f(x))} = {P_d((f(x))\pi)}$ for all $f(x)\in \mathcal{M}_c$, where $(f(x))\pi$ denotes the image of $f(x)$ under the map $\pi$. In this case, we write $c\stackrel{M}\sim d$.

$(2)$  \emph{$D$-equivalent} if there is a bijection $\pi:\mathcal{M}_c\ra\mathcal{M}_d$, such that $R[x]/(f(x)) \simeq R[x]/((f(x))\pi)$ as algebras and $\mathcal{H}_{P_c(f(x))}= \mathcal{H}_{P_d((f(x))\pi)}$ for all $f(x)\in \mathcal{M}_c$. In this case, we write $c\stackrel{D}\sim d$.

$(3)$ \emph{$AD$-equivalent} if there is a bijection $\pi:\mathcal{M}_c\ra\mathcal{M}_d$, such that $R[x]/(f(x))\simeq R[x]/((f(x))\pi)$ as algebras and either $P_c(f(x))= {P_d((f(x))\pi)}$ or $P_c(f(x))=\mathcal{J}_{P_d((f(x))\pi)}$ for all $f(x)\in \mathcal{M}_c$. In this case, we write $c\stackrel{AD}\sim d$.
\end{Def}

Clearly, $c\stackrel{M}\sim d$, $c\stackrel{D}\sim d$ and $c\stackrel{AD}\sim d$ are equivalence relations on the set of all square matrices over $R$.

Here are examples of the $D$-equivalences.
Let $R$ be a field and $J_n(\lambda)$ the $n\times n$ Jordan matrix with the eigenvalue $\lambda\in R$.

(1) We take $c=J_3(1)\oplus J_4(1)\oplus J_3(0)\oplus J_2(0)$ and $d=J_3(0)\oplus J_4(0)\oplus J_3(1)\oplus J_2(1)$. Here, $\oplus$ means forming a diagonal block matrix. In general, $m_{c\oplus d}(x)=[m_c(x),m_d(x)]$, where $[f(x),g(x)]$ stands for the least common multiple of $f(x)$ and $g(x)$ in $R[x]$. Then $m_c(x)=x^3(x-1)^4$, $\mathcal{E}_c=\{x^2, x^3, (x-1)^3, (x-1)^4\}$, $\mathcal{M}_c=\{x^3,(x-1)^4\}$, $P_c(x^3)=\{2,3\}$, $P_c((x-1)^4)=\{3,4\}$, and $m_d(x)=x^4(x-1)^3$, $\mathcal{E}_d=\{x^3,x^4,(x-1)^2,(x-1)^3\}$, $\mathcal{M}_d=\{x^4, (x-1)^3\}$, $P_d(x^4)=\{3,4\}$, $P_d((x-1)^3)=\{2,3\}$. Let $\pi: \mathcal{M}_c\to \mathcal{M}_d$ be the map: $x^3\mapsto (x-1)^3, (x-1)^4\mapsto x^4$. Then $c\stackrel{M}\sim d$. Note that $c$ and $d$ are not conjugate since they have different minimal polynomials.

(2) Let $a:=J_5(0)\oplus J_4(0)\oplus J_2(0)\in M_{11}(R)$ and $b:=J_5(0)\oplus J_3(0)\oplus J_1(0)\in M_9(R)$. Then $\mathcal{E}_a=\{x^2, x^4, x^5\}$, $\mathcal{E}_b=\{x, x^3, x^5\}$, $\mathcal{M}_a=\mathcal{M}_b=\{x^5\}$, $P_a(x^5)=\{2,4,5\}, P_b(x^5)=\{1,3,5\}$ and $\mathcal{H}_{P_a(x^5)}$ = $\{\{1,2,2\}\}$ = $ \mathcal{H}_{P_b(x^5)}.$ By definition, $a\stackrel{D}{\sim}b$, but $a\stackrel{M}{\not\sim} b$.

\subsection{Representation-finite centralizer matrix algebras\label{sect2.3}}

In this subsection we characterize representation-finite centralizer matrix algebras.

\begin{Lem}\label{iso-pr}
For $c\in M_n(R)$, the following hold true.

$(1)$ There are isomorphisms of $R$-algebras: $S_n(c,R)\simeq S_n(c^{tr},R)\simeq S_n(c,R)^{\opp}\simeq \End_{A_c}(R^n),$ where $c^{tr}$ denotes the transpose of the matrix $c$.

$(2)$ Let $\chi_c(x)$ be the characteristic polynomial of $c$. Then $S_n(c,R)=R[c]$ if and only if $\chi_c(x)=m_c(x)$.
\end{Lem}

{\it Proof.} (1) The first isomorphism follows from the fact that any matrix over a field is similar to its transpose \cite[Theorem 66, p.76]{kaplansky}, the second isomorphism is given by sending a matrix in $S_n(c^{tr},R)$ to its transpose in $S_n(c,R)^{\opp}$, and the last isomorphism follows by interpreting $c$ as a linear transformation on the $n$-dimensional $R$-space $R^n$.

(2) This is follows from Frobenius's dimension formula (see Section \ref{Introduction}). 
$\square$

\medskip
In general, $S_n(c,R)$ is neither equal to $R[c]$, nor representation-finite (see Example \ref{Bsp2.22}(2) below). But we point out when $S_n(c,R)$ is representation-finite.

\begin{Lem}\label{rep-f}
Suppose that $R$ is a perfect field, $c\in M_n(R)$ and $g(x)\in \mathcal{M}_c$. Let $b_{g(x)}:= \mbox{max}\{P_c(g(x))\cup\{3\}\}$. Then $S_n(c,R)$ is representation-finite if and only if $P_c(g(x))\subseteq \{1,b_{g(x)}-1,b_{g(x)}\}$ for all $g(x)\in \mathcal{M}_c$.
\end{Lem}

{\it Proof.} Clearly, $S_n(c,R)$ is representation-finite if and only if every block of $S_n(c,R)$ is representation-finite. The blocks of $S_n(c,R)$ are parameterized by $\mathcal{M}_c$. Let $g(x)^s\in \mathcal{M}_c$ with $g(x)\in R[x]$ an irreducible polynomial and $s\in\mathbb{N}$. Then $b_{g(x)^s}=\max\{3,s\}$ by definition.
Since $g(x)^s$ lies in $\mathcal{M}_c$, the algebra $R[x]/(g(x)^s)$ is a block of  $A_c:=R[x]/(m_c(x))$. Let $M$ be the component of the $A_c$-module $R^n$, which belongs to the block $R[x]/(g(x)^s)$, that is, $M$ is the sum of those indecomposable direct summands of $R^n$ that belong to the block $R[x]/(g(x)^s)$. Then $\End_{R[x]/(g(x)^s)}(M)$ is a block of the endomorphism algebra $\End_{A_c}(R^n)$. By Lemma \ref{add}, $\End_{R[x]/(g(x)^s)}(M)$ is Morita equivalent to $\End_{R[x]/(g(x)^s)}(\mathcal{B}(M))$. According to Lemma \ref{bijection}, $\mathcal{B}(M)\simeq \bigoplus_{t\in P_c(g(x)^s)}R[x]/(g(x)^t)$ as $A_c$-modules. Thus it follows from $S_n(c,R)\simeq \End_{A_c}(R^n)$ that each block of $S_n(c,R)$ is Morita equivalent to
$$E_{g(x)^s}:=\End_{R[x]/(g(x)^s)}\big(\bigoplus_{t\in P_c(g(x)^s)}R[x]/(g(x)^t)\big)$$for some $g(x)^s\in \mathcal{M}_c$.

Since $R$ is a perfect field, the algebraic closure  $\overline{R}$ of $R$ is a separable extension of $R$. By \cite[Theorem 3.3]{JL} which says that, for a separable extension $L/R$ of fields, a finite-dimensional $R$-algebra $\Lambda$ is representation-finite if and only if so is the $L$-algebra $L\otimes_R \Lambda$.  Hence it suffices to consider when $\overline{R}\otimes_R E_{g(x)^s}$ is representation-finite. Since $R$ is a perfect field, all irreducible factors of $m_c(x)$ are separable over $R$. Suppose
$g(x)=(x-\alpha_1)\cdots (x-\alpha_m)$, where $\alpha_1,\cdots,\alpha_m\in \bar{R}$ are pairwise distinct. Then
$$\begin{array}{ll}\overline{R}\otimes _R E_{g(x)^s} & \simeq \End_{\overline{R}\otimes_R R[x]/(g(x)^s)}\big(\overline{R}\otimes_R \bigoplus_{t\in P_c(g(x)^s)}R[x]/(g(x)^t)\big)\\ & \simeq \End_{\overline{R}[x]/(\prod^m_{i=1}(x-\alpha_i)^s)}\big(\bigoplus_{t\in P_c(g(x)^s)}\overline{R}[x]/(\prod^m_{i=1}(x-\alpha_i)^t)\big).\end{array}$$
Thus each block of $\overline{R}\otimes _R E_{g(x)^s}$ is isomorphic to $\End_{\overline{R}[x]/(x^s)}\big(\bigoplus_{t\in P_c(g(x)^s)}\overline{R}[x]/(x^t)\big). $ Now, it follows from \cite[Theorem 2.1 (i)]{YV} (see also \cite{dr}) that the endomorphism algebra $\End_{\overline{R}[x]/(x^s)}\big(\bigoplus_{t\in P_c(g(x)^s)}\overline{R}[x]/(x^t)\big) $ is representation-finite if and only if either $s\le 3$ and $P_c(g(x)^s)\subseteq \{1,2,3\}$ or $s\geq 4$ and $P_c(g(x)^s)\subseteq \{1,s-1,s\}$. This is equivalent to saying that $P_c(g(x)^s)\subseteq \{1,b_{g(x)^s}-1,b_{g(x)^s}\}$. $\square$

\medskip
As a corollary of Lemma \ref{rep-f}, we have the following.

\begin{Koro}\label{rep-p}
Let $R$ be a perfect field of characteristic $p\ge 0$, and let $\sigma\in \Sigma_n$ be a permutation of cycle type $(\lambda_1,\cdots,\lambda_s)$. Then $S_n(c_{\sigma},R)$ is representation-finite if and only if there exists a positive integer $t$ such that $\nu_p(\lambda_i)\in \{0,t\}$ for all $i\in [s]$.
\end{Koro}

{\it Proof.} Let $c:=c_{\sigma}\in M_n(R)$. If $p=0$, then $\nu_p(\lambda_i)=0$ for all $i\in [s]$. In this case, $S_n(c,R)$ is semisimple, and hence  representation-finite. Actually, let $G$ be the subgroup of $\Sigma_n$ generated by $\sigma$. Then the group algebra $R[G]$ is semisimple. Since there is a surjective homomorphism from the algebra $R[G]$ to the algebra $R[c]$ by sending $\sigma$ to $c$, we see that $R[c]$ is semisimple. Hence $ S_n(c,R)\simeq\End_{R[c]}(R^n)$ is semisimple. Thus Corollary \ref{rep-p} is true for $p=0$.

Now, we assume $p>0$.
By Lemma \ref{per}, for $g(x)\in \mathcal{M}_c$, all the integers in $P_c(g(x))$ are $p$-powers and the polynomial $(x-1)^{p^{\nu_p(\lambda_i)}}$ is an elementary divisor of $c$ for $i\in [s]$. Let $m:=\max\{\nu_p(\lambda_i)\mid i\in [s]\}$. Then $(x-1)^{p^m}\in \mathcal{M}_c$ and $P_c((x-1)^{p^m})=\{p^{\nu_p(\lambda_i)}\mid i\in [s]\}$.

Suppose that $S_n(c,R)$ is representation-finite. By Lemma \ref{rep-f}, we deduce that $P_c((x-1)^{p^m})$ does not contain two different $p$-powers $p^a>1$ and $p^b>1$ with $a\neq b$. Since $p^{\nu_p(\lambda_i)}\in P_c((x-1)^{p^m})$ for $i\in [s]$, there do not exist $\lambda_i$ and $\lambda_j$ with $i,j\in [s]$ such that $\nu_p(\lambda_i)>\nu_p(\lambda_j)\geq 1$, that is, there exists an integer $t>0$ such that $\nu_p(\lambda_i)\in \{0,t\}$ for all $i\in [s]$.

Conversely, suppose that there exists an integer $t>0$ such that $\nu_p(\lambda_i)\in \{0,t\}$ for all $i\in [s]$. Then, for $g(x)\in \mathcal{M}_c$, we deduce from Lemma \ref{per} that $P_c(g(x))\subseteq \{1,p^t\}$. Thus it follows from Lemma \ref{rep-f} that $S_n(c,R)$ is representation-finite. $\square$

Now we give nontrivial examples of representation-finite and -infinite centralizer matrix algebras.

\begin{Bsp}\label{Bsp2.22}{\rm
(1) Let $R$ be a field of characteristic $3$ and $\sigma=(123)(45)\in\Sigma_5$. Then $A:=S_5(c_{\sigma},R)$ is representation-finite by Corollary \ref{rep-p}. Now, we work out the quiver and relations for $A$. Let $f_1:=e_{11}+e_{22}+e_{33}$, $f_2:=e_{44}+e_{55}$, $f_{45}:=e_{45}+e_{54}, h_{21}:=f_{45}-f_2$ and $h_{22}:=-f_2-f_{45}$. Then $f_2=h_{21}+h_{22}$ and the set $\{f_1, h_{21}, h_{22}\}$ is a complete set of primitive orthogonal idempotents of $A$. Hence $_AA=Af_1\oplus Af_2=Af_1\oplus Ah_{21}\oplus Ah_{22}$. By calculations, we have $\dim_R(Af_1)=4$, $\dim_R(Ah_{21})=2$, $\dim_R(Ah_{22})=1$, $\dim_R(h_{22}Af_1)$ = $\dim_R(f_1Ah_{22})=0$, dim$_R(h_{21}Ah_{22}) = \dim_R(h_{22}Ah_{21})=0, \dim_R(f_1Af_1)=3$, $\dim_R(f_1Af_{2})$ $= 1$, $\dim_R(h_{22}Af_{1})=0$ and $\dim_R(h_{22}Af_{2})=1$. Let $e_3:=h_{22}, e_2:=h_{21}, e_1:=f_1, \epsilon=f_1-h_{21}, \alpha=e_{14}+e_{2,5}+e_{3,4}+e_{15}+e_{24}+e_{35}$ and $\beta=-\alpha^{tr}$. Then $A$ can be represented by the quiver with relations
$$A: \quad \xymatrix{
3\bullet & 2\bullet\ar@<2.5pt>[r]^{\beta} &\bullet 1\ar@<2.5pt>[l]^{\alpha}\ar@(ur,dr)[]^{\epsilon}
  & \qquad \alpha\beta=\epsilon^2, \; \epsilon\alpha=\beta\epsilon=\beta\alpha=0.}\vspace{-0.3cm}
$$
The Loewy structures of the indecomposable projective $A$-modules $P(i)$ are visually pictured as follows:
\[\xymatrix@C=0.1cm@R=0.2cm{
 P(1):&  & 1\ar@{-}[dl]\ar@{-}[dr] &    &  \qquad   P(2):& &   2\ar@{-}[d]     &\quad   & P(3): &\quad 3\\
 & 1 \ar@{-}[dr]    &   & 2 \ar@{-}[dl]        &  &      &  1  & & \\
 &      & 1 &               &     & & & \\
}\vspace{-0.3cm}
\]
Since $A/\soc(Ae_1)$ is representation-finite and $A$ has one more non-isomorphic indecomposable module than $A/\soc(Ae_1)$ does, $A$ is representation-finite.

(2) Let $R$ be an algebraically closed field of characteristic $2$, and let $\sigma=(1234)(56)\in \Sigma_6$ and $c:=c_{\sigma}\in M_6(R)$. Then $\mathcal{E}_c=\{(x-1)^4,(x-1)^2\}$ by Lemma \ref{per}, $m_c(x)=(x-1)^4$ and $R[c]\simeq R[x]/((x-1)^4)$.  Then the $R[c]$-module $R^6$ is isomorphic to $R[x]/((x-1)^4)\oplus R[x]/((x-1)^2)$ by $(\star)$ in Section \ref{2.3}. Hence
$$\begin{array}{ll} S_6(c,R) & \simeq \End_{R[c]}(R^6)\simeq \End_{R[x]/((x-1)^4)}\big(R[x]/((x-1)^4)\oplus R[x]/((x-1)^2)\big)\\ & \simeq \End_{R[x]/(x^4)}\big(R[x]/(x^4)\oplus R[x]/(x^2)\big).\end{array}$$
By calculations, the algebra $A:=\End_{R[x]/(x^4)}(R[x]/(x^4)\oplus R[x]/(x^2))$ can be represented by the quiver with relations:
$$\xymatrix{
	\bullet \ar@(ul,dl)[]_{\gamma} \ar@<-0.4ex>[r]_{\alpha}_(0){1}_(1){2}
	& \bullet \ar@(ur,dr)[]^{\eta} \ar@<-0.4ex>[l]_\beta}, \qquad \eta^2=\beta\alpha=0,\; \gamma^2=\alpha\beta,\; \beta\gamma=\eta\beta, \; \alpha\eta=\gamma\alpha.
$$
The Loewy structures of the indecomposable projective $A$-modules $P(1)$ and $P(2)$ can be pictured:
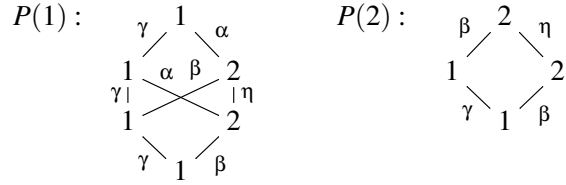
\[\xymatrix@C=0.3cm@R=0.2cm{
 P(1): &      & 1\ar@{-}[dl]_{\gamma}\ar@{-}[dr]^{\alpha} &       &  \qquad   P(2): &       &   2\ar@{-}[dl]_{\beta}\ar@{-}[dr]^-{\eta} &    \\
       & 1 \ar@{-}[d]_-{\gamma}\ar@{-}[drr]^(.3){\alpha}    &                      & 2\ar@{-}[dll]_(0.3){\beta} \ar@{-}[d]^-{\eta} &          &  1\ar@{-}[dr]_{\gamma}    &                 & 2\ar@{-}[dl]^{\beta} \\
       & 1\ar@{-}[dr]_{\gamma}    &          &   2\ar@{-}[dl]^-{\beta}            &     &  & 1 &  \\
       &     & 1 & & & &  &\\
}
\]
One can easily check that $A/\rad^2(A)$ is representation-infinite, and therefore $S_6(c,R)\simeq A $ is representation-infinite. This also follows from Corollary \ref{rep-p}.
}
\end{Bsp}

\section{Derived equivalences and homological conjectures \label{Pf}} 
This section is devoted to proving all results mentioned in the introduction.

Assume that the characteristic of $R$ is $p\geq 0$. Recall that, for $c\in M_n(R)$,  we write $A_c:=R[x]/(m_c(x))$, where  $m_c(x)$ is the minimal polynomial of $c$ over $R$. Now, let $d\in M_m(R)$, we assume the following:
$$m_c(x)=\prod^{l_c}_{i=1} f_i(x)^{n_i} \mbox{  for } n_i\ge 1 \; \mbox{ and } \; m_d(x)=\prod^{l_d}_{j=1} g_j(x)^{m_j} \; \mbox{  for } m_j\ge 1,$$   $$U_i:=R[x]/(f_i(x)^{n_i})\mbox{  for } i\in [l_c] \; \mbox{ and }\; \; V_j:=R[x]/(g_j(x)^{m_j})\mbox{ for } j\in [l_d],$$ where $f_1(x), \cdots, f_{l_c}(x)$ are pairwise distinct monic irreducible polynomials in $R[x]$, and where $g_1(x),\cdots, g_{l_d}(x)$ are pairwise distinct monic irreducible polynomials in $R[x]$. Then $U_i$ and $V_j$ are local, symmetric Nakayama $R$-algebras, and
$$A_c\simeq U_1\times U_2\times\cdots\times U_{l_c} \mbox{ and }\; A_d \simeq V_1\times V_2\times\cdots\times V_{l_d}.$$

Recall that $A_c\simeq R[c]$ and $R^n$ is viewed as an $A_c$-module. According to these blocks of $A_c$ and $A_d$, we decompose the $A_c$-module $R^n$ and the $A_d$-module $R^m$ as
$$R^n = \bigoplus^{l_c}_{i=1} M_i \; \mbox{ and } \; R^m =\bigoplus^{l_d}_{j=1} N_j,$$where $M_i$ is the sum of indecomposable direct summands of $R^n$ belonging to the block $U_i$, and where $N_j$ is the sum of indecomposable direct summands of $R^m$  belonging to the block $V_j$. Then it follows from Lemma \ref{bijection} that
$$(\dag)\quad \mathcal{B}(M_i)\simeq \bigoplus_{r\in {P_c(f_i(x)^{n_i})}} R[x]/(f_i(x)^r) \mbox{ \; and \; } \mathcal{B}(N_j)\simeq \bigoplus_{s\in {P_d(g_j(x)^{m_j})}} R[x]/(g_j(x)^s)$$as $U_i$-modules and $V_j$-modules, respectively. Since $R^n$ is a faithful $M_n(R)$-module, $R^n$ is also a faithful $R[c]$-module, and therefore $M_i$ is a faithful $U_i$-module for $i\in[l_c]$. Similarly, $N_j$ is a faithful $V_j$-module for $j\in [l_d]$.
Further, we set $$A_i:=\End_{U_i}(M_i) \; \mbox{  and } \; B_j:=\End_{V_j}(N_j)$$ for $i\in [l_c]$ and $j\in [l_d]$. Then $A_i$ and $B_j$ are indecomposable as algebras for $i\in [l_c]$ and $j\in [l_d]$. Clearly,  $A_i$ (respectively, $B_j$) is semisimple if and only if $n_i=1$ (respectively, $m_j=1$). In this case, $A_i\simeq M_k(R[x]/(f_i(x)))$ and $B_j\simeq M_t(R[x]/(g_j(x)))$ , where $k$ and $t$ are the multiplicities of $f_i(x)$ and $g_j(x)$) occurring as elementary divisors of $c$ and $d$, respectively. By Lemma \ref{iso-pr},
$$S_n(c,R)\simeq \prod^{l_c}_{i=1}{\End}_{U_i}(M_i)= \prod^{l_c}_{i=1}A_i \, \mbox{ \; and } \; S_m(d,R)\simeq \prod^{l_d}_{j=1}{\End}_{V_j}(N_j)= \prod^{l_d}_{i=1}B_j.$$

As the $R[c]$-module $R^n$ is a generator, we see that the bimodule $_{R[c]}R^n_{S_n(c,R)}$ has the double centralizer property, that is, $\End_{S_n(c,R)}(R^n_{_{S_n(c,R)}})=R[c]$.

\subsection{Characterizations of Morita and derived equivalences: Proof of Theorem \ref{main1}}
In this subsection we prove the main result, Theorem \ref{main1}.

\begin{Lem} \label{3.1} $(1)$ $\mathcal{M}_c=\{f_i(x)^{n_i}\mid  i\in [l_c]\}$.

$(2)$ If $A_i$ and $B_j$ are derived equivalent, then $U_i\simeq V_j$ and $n_i=m_j$.
\end{Lem}

{\it Proof.} (1) follows by definition. (2) is a consequence of Lemma \ref{iso}. $\square$

\begin{Lem}\label{E}
Let $c\in M_n(R)$ and $d\in M_m(R)$. Then $c\stackrel{M}\sim d$ if and only if there is an isomorphism $\varphi:R[c]\simeq R[d]$ of algebras such that $\mathcal{B}(R^n)\simeq \mathcal{B}(R^m)$, where $R^m$ is viewed as an $R[c]$-module via $\varphi$.
\end{Lem}

{\it Proof.} Suppose $c\stackrel{M}\sim d$. By definition, there is a bijection $\pi:\mathcal{M}_c\to\mathcal{M}_d$ such that, for any  $f(x)^{n_i}\in \mathcal{M}_c,$ the isomorphism $R[x]/(f(x)^{n_i})\simeq R[x]/((f(x)^{n_i})\pi)$ as algebras and ${P_c(f(x)^{n_i})}= {P_d((f(x)^{n_i})\pi)}$. Then $l_c=l_d$. It follows from
$$R[c]\simeq \prod_{f(x)^{n_i}\in \mathcal{M}_c} R[x]/(f(x)^{n_j}) \mbox{  and }\; R[d]\simeq \prod_{g(x)^{m_j}\in \mathcal{M}_d} R[x]/(g(x)^{m_j})$$that there is an isomorphism $\varphi: R[c]\simeq R[d].$ After reordering the factors in the above products, we may assume that $(f_i(x)^{n_i})\pi=g_i(x)^{m_i}$ for $i\in [l_c].$ Then the condition ${P_c(f(x)^{n_i})}= {P_d((f(x)^{n_i})\pi)},$ together with $(\dag)$, implies that $\mathcal{B}(M_i)\simeq \mathcal{B}(N_i)$  for $i\in [l_c]$. Here, $N_i$ is viewed as an $R[c]$-module via $\varphi$. Hence $\mathcal{B}(R^n)\simeq \mathcal{B}(R^m)$, where $R^m$ is viewed as an $R[c]$-module via $\varphi$.

Conversely, suppose that there is an isomorphism $\varphi:R[c]\simeq R[d]$ such that $\mathcal{B}(R^n)\simeq \mathcal{B}(R^m)$ when $R^m$ is regarded as an $R[c]$-module via $\varphi$. Then $l_c=l_d$. We may assume that $\varphi$ restricts to an isomorphism $\varphi_i: U_i\simeq V_i$, that is, $R[x]/(f_i(x)^{n_i})\simeq R[x]/(g_i(x)^{m_i}) $ for $i\in [l_c].$ This implies $n_i=m_i$ for $i\in [l_c].$ Then the condition
$\mathcal{B}(R^n)\simeq \mathcal{B}(R^m)$ implies that $\mathcal{B}(M_i)\simeq \mathcal{B}(N_i)$ for $i\in [l_c].$ Due to $(\dag)$, we have ${P_c(f_i(x)^{n_i})}= {P_d(g_i(x)^{m_i})}$ for $i\in [l_c]$.  Now we define a map $\pi: \mathcal{M}_c\to \mathcal{M}_d$ by $f_i(x)^{n_i}\mapsto g_i(x)^{m_i}$ for $i\in [l_c].$ Then $\pi$ defines an $M$-equivalence $c\stackrel{M}\sim d$. $\square$

\medskip
{\bf Proof of Theorem \ref{main1}}. Recall that $$S_n(c,R)\simeq \prod^{l_c}_{i=1}{\End}_{U_i}(M_i)= \prod^{l_c}_{i=1}A_i \, \mbox{ and } \; S_m(d,R)\simeq \prod^{l_d}_{j=1}{\End}_{V_j}(N_j)= \prod^{l_d}_{i=1}B_j.$$ If $S_n(c,R)$ and $S_m(d,R)$ are Morita (or derived, or almost $\nu$-stable derived) equivalent, then they have the same number of blocks, that is, $l_c=l_d$. Further, we may assume that $A_i$ and $B_i$ are Morita (or derived, or almost $\nu$-stable derived)) equivalent and that $F_i$ is such an equivalence for $i\in [l_c]$. As $U_i$ and $V_i$ are local Nakayama algebras for $i\in [l_c]$, it follows from Lemma \ref{3.1}(2) that there is an algebra isomorphism $\varphi_i: U_i\simeq V_i$ with $n_i=m_i$ for $i\in [l_c].$ This implies that $A_c$ and $A_d$ are isomorphic via all $\varphi_i$.

By Lemma \ref{add},  $S_n(c,R)=\End_{R[c]}(R^n)$ is Morita equivalent to $\End_{R[c]}(\mathcal{B}(R^n))$. Similarly, $S_m(d,R)=\End_{R[d]}(R^m)$ is Morita equivalent to $\End_{R[d]}(\mathcal{B}(R^m))$. 

(1) Suppose $c\stackrel{M}\sim d$. Then it follows from Lemma \ref{E} that $S_n(c,R)$ and $S_m(d,R)$ are Morita equivalent. Conversely, suppose that $S_n(c,R)$ and $S_m(d,R)$ are Morita equivalent. Then it follows from Lemma \ref{add} that $\mathcal{B}(M_i)\simeq \mathcal{B}(N_i)$ if $N_i$ is regarded as a $U_i$-module via $\varphi_i$. By identifying $A_c$ and $A_d$ with $R[c]$ and $R[d]$, respectively, we have $R[c]\simeq R[d]$ and $R^m$ can be viewed as an $R[c]$-module. Thus $\mathcal{B}(R^n)\simeq \mathcal{B}(R^m)$. By Lemma \ref{E}, we have $c\stackrel{M}\sim d.$

(2) Suppose $c\stackrel{D}\sim d$. By the definition of $D$-equivalences, $A_c\simeq A_d$ as algebras and there is a map $\pi:\mathcal{M}_c\to \mathcal{M}_d$ such that $\mathcal{H}_{P_c(f_i(x)^{n_i})}= \mathcal{H}_{P_d((f_i(x)^{n_i})\pi)}$ for $f_i(x)^{n_i}\in \mathcal{M}_c.$ Without loss of generality, we assume $(f_i(x)^{n_i})\pi=g_i(x)^{m_i}$ for $i\in [l_c].$ Then $R[x]/(f_i(x)^{n_i})\simeq R[x]/(g_i(x)^{m_i})$ as algebras and $\mathcal{H}_{P_c(f_i(x)^{n_i})}= \mathcal{H}_{P_d(g_i(x)^{m_i})}$ for $i\in [l_c]$. It follows from $(\dag)$ and Remark \ref{mult-eq} that $\End_{U_i}(\mathcal{B}(M_i))$ and $\End_{V_i}(\mathcal{B}(N_i))$ are derived equivalent. Thanks to Lemma \ref{add},  $A_i$ and $B_i$ are also derived equivalent for $i\in [l_c]$. This implies that $S_n(c,R)$ and $S_m(d,R)$ are derived equivalent.

Conversely, suppose that $S_n(c,R)$ and $S_m(d,R)$ are derived equivalent. Without loss of generality, we assume that $A_i$ and $B_i$ are derived equivalent for $i\in [l_c]$. Then, by Lemma \ref{3.1}(2), there is an isomorphism $\varphi_i: U_i\simeq V_i$ of algebras such that $U_i/\rad(U_i)\simeq V_i/\rad(V_i)$, that is, $R[x]/(f_i(x))\simeq R[x]/(g_i(x))$ for $i\in [l_c]$. Let $K_i$ be a splitting field for $f_i(x)g_i(x).$ Then $K_i\otimes_R A_i$ and $K_i\otimes_RB_i$ are derived equivalent since tensor products preserve derived equivalences (see \cite[Theorem 2.1]{JR2}).

For the irreducible polynomial $f_i(x)\in R[x]$, there is a separable irreducible polynomial $u_i(x)\in R[x]$ and an integer $s_i\in \mathbb{N}$ such that $f_i(x)=u_i(x^{p^{s_i}})$ (see, for instance, \cite[Corollary 19.9]{IS}). Here, for $p=0$, we understand $p^{s_i}=1$. Similarly, there is a separable irreducible polynomial $v_i(x)$ and an integer $t_i\in \mathbb{N}$ such that $g_i(x)=v_i(x^{p^{t_i}})$. It follows from $K_i\otimes_R \big(R[x]/(f_i(x))\big)\simeq K_i\otimes_R \big(R[x]/(g_i(x))\big)$ that $s_i=t_i$ and that $u_i(x)$ and $v_i(x)$ have the same number of roots. Therefore $f_i(x), g_i(x),u_i(x)$ and $v_i(x)$ have the same number of distinct roots in $K_i$. Let $w_i$ be the number of roots of $u_i(x)$ in $K_i$. Suppose that $\alpha_{i1}, \alpha_{i2},\cdots, \alpha_{i w_i}$ are the distinct roots of $f_i(x)$ in $K_i$ and that $\beta_{i1}, \beta_{i2}, \cdots, \beta_{i w_i}$ are the distinct roots of $g_i(x)$ in $K_i$. Then $K_i\otimes_R U_i$ $ = K_i\otimes_R \big(R[x]/(f_i(x)^{n_i})\big)\simeq \prod^{w_i}_{q=1} K_i[x]/((x-\alpha_{iq})^{n_i\cdot p^{s_i}})$. Similarly, $K_i\otimes_R V_i=K_i\otimes_R \big(R[x]/(g_i(x)^{m_i})\big)\simeq \prod^{w_i}_{q=1} K_i[x]/((x-\beta_{iq})^{m_i\cdot p^{s_i}})$.

Now, we shall show the equality $\mathcal{H}_{P_c(f_i(x)^{n_i})}= \mathcal{H}_{P_d(g_i(x)^{m_i})}$. Indeed, given a $U_i$-module $R[x]/(f_i(x)^r)$, there is the following isomorphism of $\prod^{w_i}_{q=1} K_i[x]/\big((x-\alpha_{iq})^{n_ip^{s_i}}\big)$-modules:
$$K_i\otimes_R \big(R[x]/(f_i(x)^r)\big)\simeq \bigoplus^{w_i}_{q=1} K_i[x]/\big((x-\alpha_{iq})^{rp^{s_i}}\big).$$
Note that $|P_c(f_i(x)^{n_i})|$ equals the number of non-isomorphic indecomposable direct summands of $M_i$. Since $\Hom_{U_i}(M_i,-): \add(M_i)\to A_i\mbox{-proj}$ is an equivalence, we see that $|P_c(f_i(x)^{n_i})|$ equals the number of indecomposable projective $A_i$-modules, hence the number of simple $A_i$-modules. Similarly, $|P_d(g_i(x)^{n_i})|$ is equal to the number of simple $B_i$-modules. Since derived equivalent algebras have the same number of simple modules, we get $|P_c(f_i(x)^{n_i})|=|P_d(g_i(x)^{n_i})|$. Put $h_i :=|P_c(f_i(x)^{n_i})|$. For  $h_i=1$, we have $\mathcal{H}_{P_c(f_i(x)^{n_i})}= \mathcal{H}_{P_d(g_i(x)^{m_i})}$. So we may assume that $h_i\geq 2$ and $P_c(f_i(x)^{n_i})=\{u_{i1},\cdots,u_{ih_i}\}$ with $u_{i1}> \cdots > u_{ih_i}$. Since $A_i=\End_{U_i}(M_i)$ is Morita equivalent to $\End_{U_i}(\mathcal{B}(M_i))$, the algebra $K_i\otimes_R A_i$ is Morita equivalent to  the algebra $K_i\otimes_R \End_{U_i}(\mathcal{B}(M_i))\simeq  \End_{K_i\otimes_R U_i}(K_i\otimes_R \mathcal{B}(M_i))$. As $\mathcal{B}(M_i)\simeq \bigoplus_{k=1}^{h_i} R[x]/(f_i(x)^{u_{ik}})$ as $U_i$-modules, there is the following isomorphism of $\prod^{w_i}_{q=1} K_i[x]/\big((x-\alpha_{iq})^{n_ip^{s_i}}\big)$-modules:
$$K_i\otimes_R \mathcal{B}(M_i)\simeq \bigoplus^{w_i}_{q=1}\bigoplus_{k=1}^{h_i} K_i[x]/((x-\alpha_{iq})^{u_{ik}p^{s_i}}).$$
For $q\in [w_i]$, set $E_{c,i,q}:=\End_{K_i[x]/((x-\alpha_{iq})^{n_ip^{s_i}})}\big(\bigoplus_{k=1}^{h_i} K_i[x]/((x-\alpha_{iq})^{u_{ik}p^{s_i}})\big).$ Then $\End_{K_i\otimes_R U_i}(K_i\otimes_R \mathcal{B}(M_i))\simeq \prod^{w_i}_{q=1}E_{c,i,q}$ and  $E_{c,i,q}$ is a block of $\End_{K_i\otimes_R U_i}(K_i\otimes_R \mathcal{B}(M_i))$, which is isomorphic to $E_{c,i,q'}$ for all $q'\in [w_i]$. It follows that each block of $K_i\otimes_R A_i$ is Morita equivalent to $E_{c,i,q}$ for some $q\in [w_i]$.
Similarly, we write $P_d(g_i(x)^{n_i})=\{v_{i1},\cdots,v_{ih_i}\}$ with $v_{i1}> \cdots > v_{ih_i}$, and have the following isomorphism of $\prod^{w_i}_{q=1} K_i[x]/\big((x-\beta_{iq})^{n_ip^{s_i}}\big)$-modules
$$K_i\otimes_R \mathcal{B}(N_i)\simeq \bigoplus^{w_i}_{q=1} \bigoplus_{k=1}^{h_i} K_i[x]/\big((x-\beta_{iq})^{v_{ik}p^{s_i}}\big).$$
For $q'\in [w_i]$, set $E_{d,i,q'}:=\End_{K_i[x]/((x-\beta_{iq'})^{n_ip^{s_i}})}\big(\bigoplus_{k=1}^{h_i} K_i[x]/((x-\beta_{iq'})^{v_{ik}p^{s_i}})\big).$
Then $\End_{K_i\otimes_R U_i}(K_i\otimes_R \mathcal{B}(N_i))\simeq \prod^{w_i}_{q'=1}E_{d,i,q'}$ and  $E_{d,i,q'}$ is a block of $\End_{K_i\otimes_R V_i}(K_i\otimes_R \mathcal{B}(N_i))$, which is isomorphic to $E_{d,i,q''}$ for all $q''\in [w_i]$. It follows that each block of $K_i\otimes_R B_i$ is Morita equivalent to $E_{d,i,q'}$ for some $q'\in [w_i]$.

Since $K_i\otimes_R A_i$ and $K_i\otimes_RB_i$ are derived equivalent and since derived equivalences preserve blocks, we see that $E_{c,i,q}$ and $E_{d,i,q'}$ are derived equivalent.
Note that $u_{i1}=n_i=m_i=v_{i1}$, and we have the following isomorphisms of algebras:
$$E_{c,i,q}\simeq \End_{K_i[x]/(x^{n_ip^{s_i}})}\big(\bigoplus_{k=1}^{h_i} K_i[x]/(x^{u_{ik}p^{s_i}})\big)\mbox{ \; and \; }  E_{d,i,q'}\simeq\End_{K_i[x]/(x^{n_ip^{s_i}})}\big(\bigoplus_{k=1}^{h_i} K_i[x]/(x^{v_{ik}p^{s_i}})\big).$$
Then the Cartan matrices of  $E_{c,i,q}$  and $E_{d,i,q'}$ (as $K_i$-algebras) are the $h_i\times h_i$ matrices
$$H_i:=p^{s_i}\sum^{h_i}_{k=1} (\sum^k_{l=1} u_{ik}(e_{kl}+e_{lk})-u_{ik} e_{kk}) \mbox{ and }  J_i:=p^{s_i}\sum^{h_i}_{k=1} (\sum^k_{l=1} v_{ik}(e_{kl}+e_{lk})-v_{ik} e_{kk}),$$respectively.
Since the $K_i$-algebras $E_{c,i,q}$  and $E_{d,i,q'}$ are derived equivalent and since the Cartan matrices of derived equivalent, split algebras are congruent by an invertible matrix with integral entries (see \cite[Chapter 6, Proposition 6.8.9]{Z}), there exists an invertible matrix $\Phi_i\in M_{h_i}(\mathbb{Z})$ such that $\Phi^{tr}_i H_i\Phi_i=J_i.$ Now, applying Lemma \ref{mat} to the numbers $u_{i1}>\cdots>u_{ih_i}$ and $v_{i1}>\cdots>v_{ih_i}$ as well as to the matrices $H_i$ and $J_i$, we have $\mathcal{H}_{P_c(f_i(x)^{n_i})}= \mathcal{H}_{P_d(g_i(x)^{m_i})}$ as multisets. Thus we can define a map $\pi: \mathcal{M}_c\to \mathcal{M}_d$, $f_i(x)^{n_i}\mapsto g_i(x)^{m_i}$ for $i\in [l_c].$ Then $\pi$ gives rise to a $D$-equivalence $c\stackrel{D}\sim d$.

(3) Suppose $c\stackrel{AD}\sim d.$ Then there exists a bijection $\pi: \mathcal{M}_c\to \mathcal{M}_d$, $f_i(x)^{n_i}\mapsto g_i(x)^{m_i}$ such that $\varphi_i: U_i\simeq V_i$ as algebras and  either ${P_c(f_i(x)^{n_i})}= {P_d(g_i(x)^{m_i})}$ or ${P_c(f_i(x)^{n_i})} = \mathcal{J}_{P_d(g_i(x)^{m_i})}$ for $i\in [l_c]$. By $(\dag)$ and the condition ${P_c(f_i(x)^{n_i})}= {P_d(g_i(x)^{m_i})}$ or ${P_c(f_i(x)^{n_i})} = \mathcal{J}_{P_d(g_i(x)^{m_i})}$, we have either $\mathcal{B}(M_i)_{\mathscr{P}}\simeq \mathcal{B}(N_i)_{\mathscr{P}}$ or $\mathcal{B}(M_i)_{\mathscr{P}}\simeq\Omega_{V_i}(\mathcal{B}(N_i)_{\mathscr{P}})$ as $U_i$-modules. Note that $M_i$ (respectively, $N_i$) is a faithful $U_i$-module (respectively, $V_i$-module) which contains the regular module $U_i$ (respectively, $V_i$) as a direct summand. It follows from Lemma \ref{add} that $A_i:=\End_{U_i}(M_i)$ is Morita equivalent to $\End_{U_i}\big(U_i\oplus \mathcal{B}(M_i)_{\mathscr{P}}\big)$ and that $B_i:=\End_{V_i}(N_i)$ is Morita equivalent to $\End_{V_i}\big(V_i\oplus \mathcal{B}(N_i)_{\mathscr{P}}\big)$. If $\mathcal{B}(M_i)_{\mathscr{P}}\simeq \mathcal{B}(N_i)_{\mathscr{P}}$, then $A_i$ and $B_i$ are Morita equivalent. If $\mathcal{B}(M_i)_{\mathscr{P}}\simeq\Omega_{V_i}(\mathcal{B}(N_i)_{\mathscr{P}})$, then $A_i$ and $B_i$ are almost $\nu$-stable derived equivalent by Lemma \ref{alm}. Hence, in any case, $A_i$ and $B_i$ are always almost $\nu$-stable derived equivalent, and therefore $S_n(c,R)$ and $S_m(d,R)$ are almost $\nu$-stable derived equivalent.

Conversely, suppose that $S_n(c,R)$ and $S_m(d,R)$ are almost $\nu$-stable derived equivalent. Thanks to Lemma \ref{almst}, the almost $\nu$-stable derived equivalence $F_i$ induces a stable equivalence, say $\overline{F_i}$, between $U_i$ and $V_i$, such that $\overline{F_i}(\mathcal{B}(M_i)_{\mathscr{P}})\simeq \mathcal{B}(N_i)_{\mathscr{P}}$ and $m_i=n_i$ for $i\in [l_c]$. By identifying the algebra $V_i$ with the algebra $U_i$ via $\phi_i$, we see that $\overline{F_i}$ is a stable equivalence from $U_i$ to itself. Now, according to Lemma \ref{self}, we deduce either $\mathcal{B}(M_i)_{\mathscr{P}}\simeq \mathcal{B}(N_i)_{\mathscr{P}}$ or $\mathcal{B}(M_i)_{\mathscr{P}}\simeq \Omega_{V_i}\big(\mathcal{B}(N_i)_{\mathscr{P}}\big)$ as $U_i$-modules, where $N_i$ is viewed as a $U_i$-module via $\varphi_i$. For $i\in [l_c]$, it follows from $(\dag)$ that $\mathcal{B}(M_i)_{\mathscr{P}}\simeq \mathcal{B}(N_i)_{\mathscr{P}}$ is equivalent to the condition ${P_c(f_i(x)^{n_i})}= {P_d(g_i(x)^{m_i})}$. Similarly, for $i\in [l_c]$, $\mathcal{B}(M_i)_{\mathscr{P}}\simeq \Omega_{V_i}(\mathcal{B}(N_i)_{\mathscr{P}}$ is equivalent to the condition ${P_c(f_i(x)^{n_i})}=\mathcal{J}_{P_d(g_i(x)^{m_i})}.$
Now we define a map $\pi: \mathcal{M}_c\to \mathcal{M}_d$ by $f_i(x)^{n_i}\mapsto g_i(x)^{m_i}$ for $i\in [l_c].$ By Definition \ref{newequrel}(3), $\pi$ defines an $AD$-equivalence between $c$ and $d$. $\square$

\medskip
As a corollary of Theorem \ref{main1}, we consider nilpotent matrices. For a nilpotent matrix $c\in M_n(R)$, the Jordan canonical form $c_0$ of $c$ is unique up to the ordering of its Jordan blocks. Further, $c_0$ has a Jordan block of size $t$ if and only if ${\rm rank}(c^{t+1})+{\rm rank}(c^{t-1})-2{\rm rank}(c^{t})>0$. We set $I_c:=\{t\ge 1\mid c_0~\mbox{ has a Jordan block of size t}\}.$ Note that $\mathcal{M}_c$ consists of only one polynomial of the form $x^r$ with $r$ being the maximal number in $I_c$. Thus $I_c = P_c(x^r).$

\begin{Koro}
Let $c\in M_n(R)$ be a nilpotent matrix and $d\in M_m(R)$. Then $S_n(c,R)$ and $S_m(d,R)$ are derived equivalent if and only if $d = \lambda I_m+b$ with $\lambda\in R$ and $b$ being a nilpotent matrix such that $\mathcal{H}_{I_b}=\mathcal{H}_{I_c}.$
\end{Koro}

{\it Proof.} Sufficiency. Suppose that $b\in M_m(R)$ is a nilpotent matrix and $d=\lambda I_m +b$ with $\lambda \in R$. Then $S_m(d,R)=S_m(b,R).$ Let $x^s$ be the unique polynomial in  $\mathcal{M}_b$. Furthermore, the condition $\mathcal{H}_{I_b}=\mathcal{H}_{I_c}$ implies $\mathcal{H}_{P_c(x^r)}=\mathcal{H}_{P_b(x^s)}.$ It then follows from Theorem \ref{main1} that $S_n(c,R)$ and $S_m(b,R)$ are derived equivalent.

Necessity. Suppose that $S_n(c,R)$ and $S_m(d,R)$ are derived equivalent with $c$ being nilpotent. Then  $\mathcal{M}_c=\{x^r\}$. It follows from Theorem \ref{main1} that $\mathcal{M}_d=\{h(x)^s\}$ and $R[x]/(x^r)\simeq R[x]/(h(x)^s)$ as algebras, where $h(x)$ is an irreducible monic polynomial in $R[x]$ and $s\in \mathbb{N}$. Thus $r=s$ and $h(x)=x-\lambda$ for some $\lambda\in R$. Set $b:=\lambda I_m- d$. Then $m_{b}(x)=x^s$, that is, $b$ is a nilpotent matrix. Clearly, $P_d(h(x)^s)=P_b(x^s).$ Therefore $\mathcal{H}_{P_c(x^r)}=\mathcal{H}_{P_d(h(x)^s)}=\mathcal{H}_{P_b(x^s)},$ that is, $\mathcal{H}_{I_c}=\mathcal{H}_{I_b}.$ $\square$

\medskip
Instead of $R$ being a field, we can prove the following for noetherian domains.

\begin{Rem}\label{rmk3.3}
Suppose that $R$ is a noetherian domain, $c\in M_n(R)$ and $d\in M_m(R)$. If $S_n(c,R)$ and $S_m(d,R)$ are derived equivalent, then $c\stackrel{D}\sim d$ as matrices over the fraction field of $R$.
\end{Rem}

{\it Proof.} Assume that $R$ is a noetherian domain with the fraction field $K$. Then it follows from $S_n(c,R)\subseteq M_n(R)$ that $S_n(c,R)$ is a finitely generated $R$-algebra. Thus $S_n(c,R)$ is a noetherian algebra and $S_n(c,R)\modcat$ is an abelian category, and therefore $\Db{S_n(c,R)}$ is well defined by our convention.

Regarding $K$ as an $R$-algebra, we have the isomorphism of $K$-algebras
\begin{equation}
\begin{aligned}
\varphi: K\otimes_R M_n(R)&\lra M_n(K), \; \;
\sum^s_{i=1} a_i\otimes b_i&\mapsto \, \sum^s_{i=1} (a_iI_n) b_i
\end{aligned}
\nonumber
\end{equation}
where $I_n$ is the identity matrix in $M_n(K).$ Further, $K$ is a flat $R$-module and there is the commutative diagram of $K$-algebras
$$\xymatrix{
K\otimes_RS_n(c,R) \ar[r]^-{\mu}\ar@{^{(}->}[d] & S_n(c,K)\ar@{^{(}->}[d]\\
K\otimes_R M_n(R) \ar[r]_-{\sim}^-{\varphi} & M_n(K)\\
}$$where $\mu$ is the restriction of $\varphi$. Remark that $\Img(\mu)$ belongs to $S_n(c,K)$. Since $K$ is the fraction field of $R$, we can find an element $0\ne r\in R$ for each matrix $a\in M_n(K)$ such that $ra\in M_n(R)$. This implies that $\mu$ is surjective, and therefore an isomorphism. Thus $K\otimes_R S_n(c,R)\simeq  S_n(c,K)$ as $K$-algebras.

Suppose that the $R$-algebras $S_n(c,R)$ and $S_m(d,R)$ are derived equivalent. Then there is a tilting complex $T$ over $S_n(c,R)$ such that $\End_{\Db{S_n(c,R)}}(T)\simeq S_m(d,R)$ as $R$-algebras. Since $K$ is a flat $R$-module, $\Tor^R_i(S_n(c,R),K)=0$ and $\Tor^R_i(S_m(d,R),K)=0$ for all $i\geq 1$. It then follows from \cite[Theorem 2.1]{JR2} that $K\otimes_R T$ is a tilting complex over $K\otimes_R S_n(c,R)$ with $\End_{\Db{K\otimes_R S_n(c,R)}}(K\otimes_R T)\simeq K\otimes_R S_m(d,R)$ as $K$-algebras.
Thus the $K$-algebras $S_n(c,K)$ and $S_m(d,K)$ are derived equivalent. By Theorem \ref{main1}, the equivalence $c\stackrel{D}{\sim}d$ holds as matrices over $K$.
$\square$

It is not known whether the converse of Remark \ref{rmk3.3} is true.

\subsection{Homological conjectures: Proof of Theorem \ref{fdc}\label{sect4.2}}
In this subsection, we prove that the Nakayama and finitistic dimension conjectures are true for centralizer matrix algebras.

Let $\Lambda$ be an Artin algebra, and let
$0\ra {}_{\Lambda}\Lambda \ra I_0 \ra I_1 \ra \cdots \ra I_t \ra \cdots$ be a minimal injective resolution of $_{\Lambda}\Lambda$.

\begin{Def}
$(1)$ The \emph{dominant dimension} of $\Lambda$, denoted $\dm(\Lambda)$, is the maximal $t\in \mathbb{N}$ (or $\infty$) such  that all the terms $I_0, I_1,\cdots, I_{t-1}$ in the minimal injective resolution of $_{\Lambda}\Lambda$ are projective.

$(2)$ The \emph{finitistic dimension} of $\Lambda$, denoted $\fd(\Lambda)$, is the supremum of projective dimensions of all $\Lambda$-modules $M\in \Lambda\modcat$ with finite projective dimension.
\end{Def}

Related to the two homological dimensions, there are two not yet solved major conjectures, called the \emph{Nakayama conjecture} (see \cite{Nakayama}) and the \emph{finitistic dimension conjecture} (see \cite{bass}).

\medskip
Nakayama Conjecture (NC) : An Artin algebra of infinite dominant dimension is self-injective.

Finitistic Dimension Conjecture (FDC): For any Artin algebra  $\Lambda$, $\fd(\Lambda)<\infty$.

\medskip
As is known, the validity of (FDC) for $\Lambda$ implies the validity of (NC) for $\Lambda$.
Both conjectures are open to date (see \cite[Conjectures, p.409]{ARS}). Only a few cases are verified. In the following, we will show that (FDC) holds true for all centralizer matrix algebras over fields.

\begin{Lem}{\rm \cite{it}} \label{eAe} If an Artin algebra $\Lambda$ has global dimension at most $3$, then $\fd(e\Lambda e)<\infty$ for any idempotent $e\in \Lambda$.
\end{Lem}

For a representation-finite Artin algebra $\Lambda$, let $\{X_1, \cdots, X_s\}$ be a complete set of representatives of isomorphism classes of indecomposable $\Lambda$-modules, the \emph{Auslander algebra} of $\Lambda$ is defined to be the endomorphism algebra of the $\Lambda$-module $\bigoplus_{i=1}^{s}X_i$. It is known that $\gd(\Lambda)\le 2$.

\begin{Koro} \label{rep-f-a}Let $\Lambda$ be a representation-finite Artin algebra and $A$ the Auslander algebra of $\Lambda$. Then $\fd(eAe)$ $<\infty$ for every idempotent $e\in A$. In particular, if $X\in \Lambda\modcat$, then $\fd(\End_{\Lambda}(X))<\infty$.
\end{Koro}

{\it Proof.} The first statement follows from Lemma \ref{eAe} since $\gd(A)\le 2$.  For the second statement, we may assume that $X=X^{s_1}_1\oplus\cdots \oplus X^{s_t}_t$ with $t\le s$ and  integers $s_j\ge 1$. Let $e_i\in A$ be the canonical projection from $\bigoplus_{i=1}^{s}X_i$ onto $X_i$ for $1\le i\le s$. Then $\{e_1,\cdots,e_s\}$ is a complete set of pairwise orthogonal primitive idempotent elements of $A$. Clearly, $\End_{\Lambda}(X)$ is Morita equivalent to $\End_{\Lambda}(X_1\oplus\cdots\oplus X_t)$ which is isomorphic to the algebra $(e_1+\cdots +e_t)A(e_1+\cdots +e_t)$. Again by Lemma \ref{eAe}, we get the second statement. $\square$

\medskip
Let $M$ be a generator-cogenerator for $\Lambda$-mod. The \emph{rigidity dimension} $\rd(M)$ of $M$ is defined by
$$\rd(M):=\operatorname {sup}\{n\in \mathbb{N}\mid {\Ext}^i_{\Lambda}(M,M)=0, \forall\ 1\leq i\leq n\}.$$
If no such $n$ exists, we define $\rd(M)=0$. By \cite[Lemma 3]{BJ}, $\dm(\End_\Lambda(M))=\rd(M)+2$.

The following lemma describes the dominant dimensions of centralizer matrix algebras and shows that the Nakayama conjecture holds true for centralizer matrix algebras.

Recall that, for $c\in M_n(R)$, we have a block decomposition of $S_n(c,R)$: $S_n(c,R)= \prod^{l_c}_{i=1} A_i$ with $A_i:=\End_{U_i}(M_i)$ (see the beginning of Section \ref{Pf}).

\begin{Lem}\label{lem3.1}
$(1)$ $\dm(A_i)\in \{2,\infty\}$. Particularly, $\dm(S_n(c,R))\in \{2,\infty\}$.

$(2)$ $\dm(A_i)=\infty$ if and only if $A_i$ is a symmetric, Nakayama algebra if and only if $P_c(f_i(x)^{n_i})$ is a singleton set. Thus $\dm(S_n(c,R))=\infty$ if and only if $S_n(c,R)$ is a symmetric, Nakayama algebra if and only if $P_c(f_i(x)^{n_i})$ is a singleton set for all $i\in [l_c].$
\end{Lem}

{\it Proof.} If $\Lambda$ is an Artin algebra and $L\in \Lambda\modcat$, then it follows from the Auslander-Reiten formula $D\Ext^1_{\Lambda}(L,L)\simeq \ol{\Hom}_{\Lambda}(L,\tau L)$ that $\Ext^1_{\Lambda}(L,L)\neq 0$ if $\tau L\simeq L$, where $D$ is the usual duality of an Artin algebra, $\tau := D{\rm Tr}$ denotes the Auslander--Reiten translation, and $\ol{\Hom}_{\Lambda}(X,Y)$ denotes the quotient of $\Hom_{\Lambda}(X,Y)$ modulo all homomorphisms from $X$ to $Y$ that factorize through injective $\Lambda$-modules.

Let $i\in [l_c]$. For the $U_i$-module $M_i$, we have $ (\tau M_i)_{\mathscr{P}}\simeq (M_i)_{\mathscr{P}}$, and therefore $\rd(M_i)=\infty$ if $M_i$ is projective, and $0$, otherwise. Since $\dm(A_i)=\dm(\End_{U_i}(M_i))=\rd(M_i)+2,$ we deduce that $\dm(A_i)\in \{2,\infty\}$ and that $\dm(A_i)=\infty$ if and only if $M_i$ is projective. By $(\dag)$, $M_i$ is projective if and only if $P_c(f_i(x)^{n_i})$ is a singleton set. Note that $A_i=\End_{U_i}(M_i)$ is Morita equivalent to $U_i$ if $_{U_i}M_i$ is projective. Thus $A_i$ is a symmetric, Nakayama algebra if $_{U_i}M_i$ is projective. Clearly, any symmetric algebra has infinite dominant dimension. Since $\dm(\Lambda\oplus\Gamma)=\min\{\dm(\Lambda),\dm(\Gamma)\}$ for Artin algebras $\Lambda$ and $\Gamma$, we have $\dm(S_n(c,R))=\min\{\dm(A_i)\mid i\in [l_c]\}\in \{2,\infty\}.$ Thus, $\dm(S_n(c,R))=\infty$ if and only if $\dm(A_i)=\infty$ for all $i\in [l_c]$ if and only if $P_c(f_i(x)^{n_i})$ is a singleton set for all $i\in [l_c]$ if and only if $S_n(c,R)$ is a symmetric, Nakayama algebra. $\square$

\medskip
{\bf Proof of Theorem \ref{fdc}}. (1) Let $R$ be a field and $c\in M_n(R)$. Since Nakayama algebras are representation-finite \cite[Lemma 2.1, p.197]{ARS}, their Auslander algebras have global dimension at most $2$. All blocks of $S_n(c,R)$ are of the form $A_i=\End_{U_i}(M_i)$, $i\in[l_c]$, where $U_i$ is a symmetric Nakayama algebra and $M_i$ is a generator for $U_i\modcat$. By Lemma \ref{rep-f-a}, $\fd(A_i)=\fd(\End_{U_i}(M_i))<\infty$ for all $i\in [l_c]$. Since $\fd(S_n(c,R))$ $= \mbox{max}\{\fd(A_i)\mid i\in [l_c]\}$, we see $\fd(S_n(c,R)) < \infty$. Since the validity of (FDC) for an Artin algebra $\Lambda$ implies the validity of (NC) for the same Artin algebra $\Lambda$. Hence (NC) holds true for $S_n(c,R)$. This also follows from Lemma \ref{lem3.1}(2).

(2) Let $c\in M_n(R)$ and $d\in M_m(R)$. Suppose that $S_n(c,R)$ and $S_m(d,R)$ are derived equivalent. By Lemma \ref{lem3.1}(1), $\dm(S_n(c,R)) \in \{2,\infty\}$. Thus, to prove that $S_n(c,R)$ and $S_m(d,R)$ have the same dominant dimension, we only need to show that $\dm(S_n(c,R)) = \infty$ if and only if $\dm(S_m(d,R))=\infty$. However, this follows from Theorem \ref{main1} about derived equivalences and Lemma \ref{lem3.1}(2) immediately. Thus $S_n(c,R)$ and $S_m(d,R)$ have the same dominant dimension.
$\square$

\subsection{Derived equivalences imply Morita equivalences: Proof of Corollary \ref{derp}}

To prove Corollary \ref{derp} , we recall a result on stable equivalences of Morita type.

Following \cite[Section 2]{hx3}, we say that a stable equivalence $\Phi: A\stmc{}\ra B\stmc$ of Morita type \emph{lifts} to a Morita equivalence if there is a Morita equivalence $F: A\modcat{}\ra B\modcat$ such that the  following diagram of functors is commutative (up to natural isomorphism)
$$\xy
(0,15)*+{\stmc{A}}="a",
(25,15)*+{\stmc{B}}="b",
(0,0)*+{A\modcat}="d",
(25,0)*+{B\modcat}="e",
{\ar^{can.} "d";"a"},
{\ar^{can.} "e";"b"},
{\ar^{\Phi}, "a";"b"},
{\ar^{F}, "d";"e"},
\endxy$$

Given an idempotent $e$ in an algebra $A$, the functor $Ae\otimes_{eAe}-: eAe\modcat \to A\modcat$ is called a \emph{Schur functor} that is fully faithful. Clearly, this Schur functor induces a functor on the stable module categories: $\stmc{eAe} \ra A\stmc$. For simplicity, the induced functor is still called a Schur functor.

The following lemma, taken essentially from \cite{hx3}, provides a way to get Morita equivalences from stable equivalences of Morita type (see Section \ref{sect2.1} for Definition).

\begin{Lem}\label{Sta-M}
Let $A$ and $B$ be algebras without nonzero semisimple direct summands such that $A/\rad(A)$ and $B/\rad(B)$ are separable, and let $e\in A$ and $f\in B$ be $\nu$-stable idempotents such that $eAe$ and $fAf$ are the Frobenius parts of $A$ and $B$, respectively. Suppose there is a stable equivalence $\Phi: A\stmc{}\ra B\stmc$ of Morita type. Then the following hold.

$(1)$ If $\Phi(S)$ is isomorphic in $B\stmc$ to a simple $B$-module for each simple $A$-module $S$, then $\Phi$ lifts to a Morita equivalence.

$(2)$ The functor $\Phi$ restricts to a stable equivalence $\Phi_1: eAe\stmc{}\ra fBf\stmc$ of Morita type such that the following diagram is commutative (up to natural isomorphism)
$$\xy
(0,15)*+{\stmc{A}}="a",
(25,15)*+{\stmc{B}}="b",
(0,0)*+{\stmc{eAe}}="d",
(25,0)*+{\stmc{fBf}}="e",
{\ar^{\lambda} "d";"a"},
{\ar^{\lambda} "e";"b"},
{\ar^{\Phi}, "a";"b"},
{\ar^{\Phi_1}, "d";"e"},
\endxy$$
where $\lambda$ stands for the Schur functor. Moreover, if $\Phi_1$ lifts to a Morita equivalence, then so does $\Phi$.
\end{Lem}

{\it Proof.} (1) is just \cite[Proposition 3.3]{hx3}. (2) The first statement follows from \cite[Theorem 4.2]{DM}, see also \cite[Section 3]{hx3}. The last statement follows from \cite[Proposition 3.5]{hx3}. $\square$

\medskip
{\bf Proof of Corollary \ref{derp}}.  Let $c\in M_n(R)$ and $d\in M_m(R)$.

(2) Assume that $c$ and $d$ are permutation matrices and that $S_n(c,R)$ and $S_m(d,R)$ are derived equivalent. Then $S_n(c,R)$ and $S_m(d,R)$ have the same number of blocks, that is, $l_c=l_d$. So we may assume that $A_i$ and $B_i$ are derived equivalent for $i\in [l_c].$ By Lemma \ref{3.1}, $U_i\simeq V_i$ and $n_i=m_i$ for $i\in [l_c].$ By Theorem \ref{main1} on Morita equivalences, it suffices to show that $P_c(f_i(x)^{n_i})=P_d(g_i(x)^{m_i})$ for $i\in [l_c].$

Actually, by Lemma \ref{per}, the integers in $P_c(f_i(x)^{n_i})$ and in $P_c(g_i(x)^{m_i})$ are $p$-powers for $i\in [l_c]$. We have seen in the proof of Theorem \ref{main1} about derived equivalences that $P_c(f_i(x)^{n_i})$ and $P_d(g_i(x)^{n_i})$ have the same cardinality. Let $t_i:=|P_c(f_i(x)^{n_i})|=|P_d(g_i(x)^{n_i})|$ for $i\in [l_c]$. If $t_i=1$ (this may happen for $p=0$), then $P_c(f_i(x)^{n_i})=\{n_i\}=\{m_i\}=P_d(g_i(x)^{m_i})$. Now, we may assume that $t_i\geq 2$ and $p>0$. Let $P_c(f_i(x)^{n_i}):=\{p^{u_1},\cdots,p^{u_{t_i}}\}$ with $u_1> \cdots > u_{t_i}$ and $P_d(g_i(x)^{m_i}):=\{p^{v_1},\cdots,p^{v_{t_i}}\}$ with $v_1> \cdots > v_{t_i}$. By Theorem \ref{main1} on derived equivalences, we get $\{p^{u_1}-p^{u_{2}},\cdots, p^{u_{t_i-1}}-p^{u_{t_i}},p^{u_{t_i}}\}= \{p^{v_1}-p^{v_2},\cdots, p^{v_{t_i-1}}-p^{v_{t_i}},p^{v_{t_i}}\}$. Notice the following basic facts:

(i) For integers $a>b>0$, the number $p^a-p^b$ is a $p$-power if and only if $p=2$ and $a=b+1$;

(ii) For integers $a>b>0$ and $s>t>0$, the equality $p^a-p^b=p^s-p^t$ holds if and only if $a=s$ and $b=t$.

By considering the cases $p=2$  and $p\geq 3$ separately, we get $u_k=v_k$ for all $k\in [t_i].$  Thus $P_c(f_i(x)^{n_i})=P_d(g_i(x)^{m_i})$ for $i\in [l_c]$. This implies that $A$ and $B$ are Morita equivalent by Theorem \ref{main1}.

(3) Suppose that the field $R$ is perfect. Then all irreducible factors of $m_c(x)$ are separable polynomials over $R$. Let $A_i:=\End_{U_i}(M_i)$ be a block in $S_n(c,R)$ and $P$ an arbitrary indecomposable projective $A_i$-module. Then $P\simeq  \Hom_{U_i}(M_i,X)$ for some indecomposable direct summand $X$ of the $U_i$-module $M_i$. Thus $\End_{A_i}(P)\simeq \End_{U_i}(X)$, and therefore $$\End_{A_i}(\top(P))=\End_{A_i}(P)/\rad(\End_{A_i}(P))\simeq \End_{U_i}(X)/\rad(\End_{U_i}(X)),$$where $\top(P)$ denotes the quotient $P/\rad(P)$ of a module $P$ by its radical. For the indecomposable $U_i$-module $X$, we have $\End_{U_i}(X)\simeq R[x]/(f_i(x)^t)$ for some positive integer $t$. Thus $\End_{A_i}(\top(P))\simeq \End_{U_i}(X)/\rad(\End_{U_i}(X))\simeq R[x]/(f_i(x))$ is separable.
Hence all the semisimple quotients of blocks of $S_n(c,R)$ are separable. Similarly, all the semisimple quotients of blocks of $S_m(d,R)$ are separable.

(a) $\Rightarrow$ (b) Suppose that $S_n(c,R)$ and $S_m(d,R)$ are almost $\nu$-stable derived equivalent. Then, by \cite[Theorem 1.1]{hx1}, there is a stable equivalence $F$ of Morita type between $S_n(c,R)$ and $S_m(d,R)$. Further, by Theorem \ref{main1}, we have $c\stackrel{AD}{\sim} d$, that is, there is a bijection $\pi$ between $\mathcal{M}_c$ and $\mathcal{M}_d$ such that $R[x]/(f(x))\simeq R[x]/((f(x))\pi)$ as algebras and either $P_c(f(x))= {P_d((f(x))\pi)}$ or $P_c(f(x))=\mathcal{J}_{P_d((f(x))\pi)}$ for all $f(x)\in \mathcal{M}_c$. Clearly, $\pi$ maps only irreducible polynomials to irreducible polynomials. Thus $\pi$ induces a bijection between $\mathcal{M}_c\setminus\mathcal{R}_c$ and $\mathcal{M}_d \setminus \mathcal{R}_d$  such that $R[x]/(f(x))\simeq R[x]/((f(x))\pi)$ as algebras for $f(x)\in \mathcal{M}_c\setminus\mathcal{R}_c$.

(b) $\Rightarrow$ (a) Suppose that $S_n(c,R)$ and $S_m(d,R)$ are stably equivalent of Morita type and there is a bijection $\pi: \mathcal{M}_c\setminus\mathcal{R}_c\ra \mathcal{M}_d\setminus\mathcal{R}_d$, such that $R[x]/(f(x))\simeq R[x]/((f(x))\pi)$ as algebras for $f(x)\in \mathcal{M}_c\setminus\mathcal{R}_c$. By Theorem \ref{main1}, it suffices to show $c\stackrel{AD}\sim d$. Note that an irreducible elementary divisor $f(x)$ in $\mathcal{M}_c\setminus\mathcal{R}_c$ corresponds to a semisimple block of $S_n(c,R)$, which is Morita equivalent to $R[x]/(f(x))$. Similarly, an irreducible elementary divisor $g(x)$ in $\mathcal{M}_d\setminus\mathcal{R}_d$ corresponds to a semisimple block of $S_m(d,R)$, which is Morita equivalent to $R[x]/(g(x))$. Thus the assumption on $\pi$ implies that the product of semisimple blocks of $S_n(c,R)$ and the product of semisimple blocks of $S_m(d,R)$ are Morita equivalent. Let $A_1,\cdots,A_s$ be the non-semisimple blocks of $S_n(c,R)$ with $A_i:=\End_{U_i}(M_i)$, and let $B_1,\cdots,B_t$ be the non-semisimple blocks of $S_m(d,R)$ with $B_j:=\End_{V_j}(N_j)$. Suppose that $F$ is a stable equivalence of Morita type between $S_n(c,R)$ and $S_m(d,R)$. Then $F$ induces a stable equivalence of Morita type between $\bigoplus_{i=1}^sA_i$ and $\bigoplus_{j=1}^tB_j$. Thus $s=t$ by \cite[Theorem 2.2]{Liu1}, and we may assume that $F$ induces a stable equivalence $F_i$ of Morita type between $A_i$ and $B_i$ for $i\in [s]$.

To show $c\stackrel{AD}\sim d$, we consider the generator $M_i$ for $U_i\modcat$.  It follows from $\nu_{A_i}\Hom_{U_i}(M_i, U_i) \simeq
\Hom_{U_i}(M_i,\nu_{U_i}U_{i})$ (see \cite[Remark 2.9 (2)]{hx3}) that the Frobenius parts of $A_i$ and $B_i$ are Morita equivalent to $U_i$ and $V_i$,
respectively. Since $A_i/\rad(A_i)$ and $B_i/\rad(B_i)$ are separable, we deduce from Lemma \ref{Sta-M}(2) that $F_i$ restricts to a stable equivalence $G_i$ of Morita type between $U_i$ and $V_i$. As $f_i(x)$ is separable and both $A_i$ and $B_i$ are non-semisimple, Corollary \ref{St-i} implies $U_i\simeq V_i$, that is, $R[x]/(f_i(x)^{n_i})\simeq R[x]/(g_i(x)^{m_i})$, and $n_i=m_i$.

Now we regard $V_i$-modules as $U_i$-modules via this isomorphism. Let $\overline{A}_i:=\End_{U_i}(U_i\oplus \mathcal{B}(M_i)_{\mathscr{P}})$,  $\overline{B}_i:=\End_{V_i}(V_i\oplus \mathcal{B}(N_i)_{\mathscr{P}})$ and $\overline{C}_i:=\End_{V_i}(V_i\oplus \Omega_{V_i}(\mathcal{B}(N_i)_{\mathscr{P}}))$, and let $e,f$ and $g$ be the $\nu$-stable idempotents of $\overline{A}_i,\overline{B}_i$ and $\overline{C}_i$, defining their Frobenius parts, respectively. Then any  two algebras from the list $\{A_i,\overline{A}_i,B_i,\overline{B}_i,\overline{C}_i\}$ are stably equivalent of Morita type (see Lemmas \ref{add} and \ref{alm}), and there is the following commutative (up to natural isomorphism) diagram by Lemma \ref{Sta-M}(2):
$$\xy
(0,15)*+{\stmc{\overline{A}_i}}="a",
(25,15)*+{\stmc{\overline{B}_i}}="b",
(50,15)*+{\stmc{\overline{C}_i}}="c",
(0,0)*+{\stmc{e\overline{A}_i e}}="d",
(25,0)*+{\stmc{f\overline{B}_i f}}="e",
(50,0)*+{\stmc{g\overline{C}_i g}}="f",
{\ar^{\lambda} "d";"a"},
{\ar^{\lambda} "e";"b"},
{\ar^{\lambda} "f";"c"},
{\ar^{\Phi}, "a";"b"},
{\ar^{\Phi_1}, "d";"e"},
{\ar^{\Psi}, "b";"c"},
{\ar^{\Psi_1}, "e";"f"},
\endxy$$
where $\lambda$ is the full embedding of stable module categories induced by the corresponding Schur functor and where $\Phi$ and $\Psi$ define stable equivalences of Morita type between $\overline{A}_i$ and $\overline{B}_i$, and between $\overline{B}_i$ and $\overline{C}_i$, respectively, while $\Phi_1$ and $\Psi_1$ are the restrictions of $\Phi$ and $\Psi$, respectively. They are again of Morita type (see Lemma \ref{Sta-M}(2)). Note that $e\overline{A}_i e\simeq U_i\simeq V_i\simeq f\overline{B}_i f\simeq g\overline{C}_i g$, and all of them are local symmetric, Nakayama algebras. Identifying $f\overline{B}_i f$ with $g\overline{C}_i g$, we can choose $\Psi$ so that $\Psi_1$ is the syzygy functor on $f\overline{B}_i f\stmc$ (see the arguments in \cite[Proposition 3.3 and Corollary 3.4]{LX3}). Let $S$  be the unique simple $e\overline{A}_i e$-module (up to isomorphism). If we identify $e\bar{A}_ie$ with $f\bar{B}_if$, then it follows from Lemma \ref{self} that either $\Phi_1(S)\simeq S$ or $\Phi_1(S)\simeq \Omega_{e\overline{A}_i e}(S)$. Thus either
$\Phi_1(S)$ or $\Psi_1\circ \Phi_1(S)$ is simple. By Lemma \ref{Sta-M}(1), either $\Phi_1$ or $\Psi_1\circ \Phi_1$ can be lifted to a Morita equivalence, and therefore
either $\Phi$ or $\Psi\circ \Phi$ can be lifted to a Morita equivalence by Lemma \ref{Sta-M}(2). It then follows from Lemma \ref{add} that either $\mathcal{B}(M_i)_{\mathscr{P}}\simeq \mathcal{B}(N_i)_{\mathscr{P}}$ or $\mathcal{B}(M_i)_{\mathscr{P}}\simeq \Omega_{V_i}(\mathcal{B}(N_i)_{\mathscr{P}})$.
By $(\dag)$, we have ${P_c(f_i(x)^{n_i})}= {P_d(g_i(x)^{m_i})}$ or ${P_c(f_i(x)^{n_i})}=\mathcal{J}_{P_d(g_i(x)^{m_i})}$. Now we define a map $\pi': \mathcal{M}_c\lra \mathcal{M}_d,$ $f_i(x)^{n_i}\mapsto g_i(x)^{m_i} \mbox{ for } f_i(x)^{n_i}\in \mathcal{R}_c, \; f(x)\mapsto (f(x))\pi \; \mbox{ for } f(x)\in \mathcal{M}_c\setminus\mathcal{R}_c.$ Then $\pi'$ defines an $AD$-equivalence of matrices $c$ and $d$, that is, $c\stackrel{AD}\sim d$. $\square$

\smallskip
Since Morita equivalences  preserve dominant, finitistic and global dimensions, we have the following.

\begin{Koro} If permutation matrices $c\in M_n(R)$ and $d\in M_m(R)$ are $D$-equivalent, then
$$\dm((S_n(c,R))=\dm((S_m(d,R)), \fd(S_n(c,R))=\fd(S_m(d,R))\mbox{ and }$$
$\gd(S_n(c,R))=\gd(S_m(d,R)),$ where $\gd(A)$ denotes the global dimension of an algebra $A$.
\end{Koro}

\subsection{Derived equivalences for permutation matrices: Proof of Corollary \ref{cor1.5}\label{sect4}}
In this subsection we discuss relations between derived equivalences of centralizer matrix algebras of permutation matrices on the one hand and derived equivalences of centralizer matrix algebras of permutation matrices of $p$-regular and $p$-singular parts on the other hand. This provides a proof of Corollary \ref{cor1.5}.

Given a prime number $p>0$ and a permutation $\sigma=\sigma_1\cdots\sigma_k\in \Sigma_n$, which is the product of disjoint cycle-permutations $\sigma_i$ of cycle type $\lambda=(\lambda_1,\cdots, \lambda_k)$ with $\lambda_i\ge 1$ for $i\in [k]$, we say that $\sigma_i$ is \emph{$p$-regular} if $p\nmid \lambda_i$, and \emph{$p$-singular} if $p\mid \lambda_i$.
The \emph{$p$-regular part} $r(\sigma)$ of $\sigma$ is the product of $p$-regular cycles of $\sigma$, and the \emph{$p$-singular part} $s(\sigma)$) of $\sigma$ is the product of $p$-singular cycles of $\sigma$. Both $r(\sigma)$ and $s(\sigma)$ are considered as elements in $\Sigma_n$, that is, $r(\sigma)$ fixes the elements involved in the $p$-singular cycles, and $s(\sigma)$ fixes the ones in $p$-regular cycles of $\sigma$. Let $c_{\sigma}:=\sum_{i=1}^ne_{i,(i)\sigma}\in M_n(R)$ be the permutation matrix of $\sigma$, where $e_{ij}$ is the matrix with $1$ in $(i,j)$-entry and $0$ in all other entries.

We start with the following corollary.

\begin{Koro}\label{one-more}
Let $R$ be a noetherian domain of characteristic $p>0$ and $\sigma\in \Sigma_n$ be of cycle type $\lambda:=(\lambda_1, \cdots, \lambda_k)$, and let $\sigma^+$ be a permutation in $\Sigma_{n+1}$ of cycle type $\lambda^+:=(\lambda_1, \cdots, \lambda_k, 1)$. Then the following are equivalent

$(a)$ $S_n(c_\sigma,R)$ and $S_{n+1}(c_{\sigma^+},R)$ are derived equivalent.

$(b)$ $S_n(c_\sigma,R)$ and $S_{n+1}(c_{\sigma^+},R)$ are Morita equivalent.

$(c)$ There exists a natural number $i\in [k]$ such that $p\nmid \lambda_i$.
\end{Koro}

{\it Proof.} Let $K$ be the fraction field of $R$ and $\mathbb{F}_p$ be the prime field of $K$. Since $c_{\sigma^+}$ is just the diagonal block matrix diag$(c_{\sigma},1)$, we have $\mathcal{E}_{c_{\sigma^+}}=\mathcal{E}_{c_\sigma}\cup \{x-1\}$ when $c_\sigma$ and $c_{\sigma^+}$ are viewed as matrices over either $K$ or $\mathbb{F}_p$.
Note that all $\lambda_i$ are exactly the orbit lengths of the cyclic group $\langle\sigma\rangle$ acting on $[n]$.

$(a)\Rightarrow (c)$ Suppose that $S_n(c_\sigma,R)$ and $S_{n+1}(c_{\sigma^+},R)$ are derived equivalent. Then it follows from Remark \ref{rmk3.3} that $S_n(c_\sigma,K)$ and $S_{n+1}(c_{\sigma^+},K)$ are derived equivalent, and hence Morita equivalent by Corollary \ref{derp}. It then follows from Theorem \ref{main1} that $c_\sigma\stackrel{M}\sim c_{\sigma^+}$ as matrices over $K$. Since $|\mathcal{E}_d|=\sum_{f(x)\in \mathcal{M}_d}|P_d(f(x))|$ for any matrix $d$, the $M$-equivalence between $c_\sigma$ and $c_{\sigma^+}$ implies that $|\mathcal{E}_{c_\sigma}|=|\mathcal{E}_{c_{\sigma^+}}|$. Now, it follows from $\mathcal{E}_{c_{\sigma^+}}=\mathcal{E}_{c_\sigma}\cup \{x-1\}$ that $x-1\in \mathcal{E}_{c_\sigma}$. But, by Lemma \ref{per}, $x-1\in \mathcal{E}_{c_\sigma}$ if and only if there is some $i\in [k]$ such that $p\nmid \lambda_i$.

$(c)\Rightarrow (b)$ Assume (c). Then there is some $i$ such that $p\nmid \lambda_i$. It follows from  $\nu_p(\lambda_i)=0$ and Lemma \ref{per} that $x-1\in \mathcal{E}_{c_\sigma}.$ Thus $\mathcal{E}_{c_{\sigma}}=\mathcal{E}_{c_{\sigma^+}}$. By Theorem \ref{main1}, $S_n(c_\sigma,\mathbb{F}_p)$ and $S_{n+1}(c_{\sigma^+},\mathbb{F}_p)$ are Morita equivalent. Therefore $R\otimes_{\mathbb{F}_p}S_n(c_\sigma,\mathbb{F}_p)$ and $R\otimes_{\mathbb{F}_p}S_{n+1}(c_{\sigma^+},\mathbb{F}_p)$ are Morita equivalent. With an argument similar to the one in Remark \ref{rmk3.3}, we obtain the isomorphisms of $R$-algebras
$$R\otimes_{\mathbb{F}_p}S_n(c_\sigma,\mathbb{F}_p)\simeq S_n(c_\sigma,R)\; \mbox{ and } \; R\otimes_{\mathbb{F}_p}S_{n+1}(c_{\sigma^+},\mathbb{F}_p)\simeq S_{n+1}(c_{\sigma^+},R).$$ Hence $S_n(c_\sigma,R)$ and $S_{n+1}(c_{\sigma^+},R)$ are Morita equivalent.

$(b)\Rightarrow (a)$ This is obvious. $\square$

\begin{Prop}\label{regular-singular}
Let $R$ be a field of characteristic $p\ge 0$, $\sigma\in \Sigma_n$ and $\tau\in \Sigma_m$. If $S_n(c_{\sigma},R)$ and $S_m(c_{\tau},R)$ are derived equivalent, then

$(1)$ $S_n(c_{r(\sigma)},R)$ and $S_m(c_{r(\tau)},R)$ are derived equivalent, and

$(2)$ $S_n(c_{s(\sigma)},R)$ and $S_m(c_{s(\tau)},R)$ are derived equivalent.
\end{Prop}

{\it Proof.} Let $\lambda=(\lambda_1,\cdots,\lambda_k)$ be the cycle type of $\sigma$. We know from Lemma \ref{per} that $m_{c_\sigma}(x)={\rm lcm}(x^{\lambda_1}-1,\cdots, x^{\lambda_k}-1)$, the least common multiple of $x^{\lambda_i}-1$, $i\in [k]$.
Recall that $\nu_p(n)$ denotes the largest non-negative integer such that $p^{\nu_p(n)}$ divides $n$, and for an irreducible factor $f(x)$ of $m_{c_{\sigma}}(x)$,  we define $$q_{f(x)}:= \mbox{max}\{\nu_{p}(\lambda_j) \mid j\in[k] \mbox{ such that } f(x)\mbox{ divides } x^{\lambda_j}-1\}.$$
According to Lemma \ref{per}, we have
$$\quad (\alpha) \quad \mathcal{E}_{c_{\sigma}} = \{ f(x)^{p^{\nu_p(\lambda_i)}}\mid i\in [k], f(x)~\mbox{is an irreducible factor of}~x^{\lambda_i}-1\}, \mbox{ and } $$ $$(\beta) \quad \mathcal{M}_{c_\sigma}=\{f(x)^{p^{q_{f(x)}}}\mid f(x)~\mbox{is an irreducible factor of }~m_{c_{\sigma}}(x)\}.\qquad \qquad$$In particular,

\smallskip
$\quad \qquad \quad \; \; ( \gamma )\quad \smallskip
\mathcal{M}_{c_\sigma}$ always contains an elementary divisor $(x-1)^{p^a}$ for some integer $a\geq 0$.

Note that $x-1\notin  \mathcal{E}_{c_{\sigma}}$ if and only if $\nu_p(\lambda_i)>0$ for all $i\in [k]$ if and only if $\sigma=s(\sigma)$ if and only if $r(\sigma)=id$, the identity permutation in $\Sigma_n$.

If $p=0$, then the statements (1) and (2) are trivially true. In the following, we assume $p>0$.

Let $\{\lambda_{j_1},\cdots, \lambda_{j_l}\}$ be the set of parts $\lambda_i$ of $\lambda$ such that $p\nmid \lambda_i$,  and let $\{\lambda_{i_1},\cdots,\lambda_{i_t}\}$ be the set of parts $\lambda_j$ of $\lambda$ such that $\nu_p(\lambda_j)>0$. We define $\ell_r(\lambda):=\sum_{i=1}^l\lambda_{j_i}$ and $\ell_s(\lambda):=\sum_{j=1}^t\lambda_{i_j}.$  Then $n=\sum_{i=1}^k\lambda_i=\ell_r(\lambda)+\ell_s(\lambda)$. The cycle type of  $r(\sigma)$ is $(\lambda_{j_1},\cdots, \lambda_{j_l}, \underbrace{1, \cdots, 1}_{n-\ell_r(\lambda)})$, and the cycle type of  $s(\sigma)$ is $(\lambda_{i_1},\cdots, \lambda_{i_t},\underbrace{1,\cdots, 1}_{n-\ell_s(\lambda)})$.

It follows from $(\alpha)$ and $(\beta)$ that
$$\begin{array}{ll}\mathcal{E}_{c_{r(\sigma)}} = \mathcal{M}_{c_{r(\sigma)}} &= \{f(x)\in R[x]\mid \exists\; a\in [l], f(x)~\mbox{is an irreducible factor of }~x^{\lambda_{j_a}}-1\} \\ &  \cup\{x-1\} \\
&= \{f(x)\in \mathcal{E}_{c_\sigma}\mid f(x)~\mbox{is irreducible }\}\cup \{x-1\},
\end{array}$$
$$\mathcal{E}_{c_{s(\sigma)}}=\begin{cases}\{u(x)\in \mathcal{E}_{c_\sigma}\mid u(x)~\mbox{is reducible in}~ R[x]\} & \mbox{ if } s(\sigma)=\sigma,\\ \{u(x)\in \mathcal{E}_{c_\sigma}\mid u(x)~\mbox{is reducible in}~ R[x]\}\cup \{x-1\} & \mbox{ if } s(\sigma)\neq \sigma.\end{cases}$$and $$\mathcal{M}_{c_{s(\sigma)}}=\begin{cases}\{g(x)\in \mathcal{M}_{c_\sigma}\mid g(x)~\mbox{ is reducible }\} & \mbox{ if } s(\sigma)\neq id,\\\{x-1\} & \mbox{ if } s(\sigma)=id.\end{cases}$$

Thus, we have the following for the power index sets.

\smallskip
$\quad (\delta)$ If $s(\sigma)\neq id$, then $P_{c_{s(\sigma)}}(h(x))=P_{c_\sigma}(h(x))\setminus\{1\}$ for $h(x)\in \mathcal{M}_{c_{s(\sigma)}}\setminus\{(x-1)^{p^a}\}$ and $P_{c_\sigma}((x-1)^{p^a})=P_{c_{s(\sigma)}}((x-1)^{p^a})$. Similar conclusions hold for $\tau\in \Sigma_m$.

\smallskip
Suppose $c_{\sigma}\stackrel{D}\sim c_{\tau}$, that is, $c_{\sigma}\stackrel{M}\sim c_{\tau}$ by Corollary \ref{derp}(2). Then there is a bijection $\pi: \mathcal{M}_{c_\sigma}\to \mathcal{M}_{c_\tau}$ such that $R[x]/(h(x))\simeq R[x]/((h(x)\pi)$ as algebras and $P_{c_\sigma}(h(x))=P_{c_\tau}((h(x))\pi)$ for $h(x)\in \mathcal{M}_{c_\sigma}$. We show that

(i) $r(\sigma)=id$ if and only if $r(\tau)=id$.

(ii) $s(\sigma)=id$ if and only if $s(\tau)=id$.

In fact, for nonnegative integers $a,b$, if $R[x]/(w(x)^a)\simeq R[x]/(z(x)^b)$ as algebras for two irreducible polynomials $w(x),z(x)\in R[x]$, then $a=b$ and $R[x]/(w(x)^i)\simeq R[x]/(z(x)^i)$ as algebras for all $i\le a$. Thus we may extend $\pi$ to a bijection between $\mathcal{E}_{c_\sigma}$ and $\mathcal{E}_{c_\tau}$ such that $R[x]/(h(x))\simeq R[x]/((h(x)\pi)$ as algebras for $h(x)\in \mathcal{E}_{c_\sigma}$. Note that $x-1\notin  \mathcal{E}_{c_{\sigma}}$ if and only if $\nu_p(\lambda_i)>0$ for $i\in [k]$ if and only if $\sigma=s(\sigma)$ if and only if $1\notin P_{c_\sigma}(h(x))$ for $h(x)\in \mathcal{M}_{c_\sigma}$. Similarly, the above observation holds for $\tau$. Thus we deduce from $P_{c_\sigma}(h(x))=P_{c_\tau}((h(x))\pi)$ for all $h(x)\in \mathcal{M}_{c_\sigma}$ that $x-1\notin  \mathcal{E}_{c_{\sigma}}$ if and only if $x-1\notin  \mathcal{E}_{c_{\tau}}$. This implies that $r(\sigma)=id$ if and only if $r(\tau)=id$.
Note that $s(\sigma)=id$ if and only if $r(\sigma)=\sigma$ if and only if  $P_{c_\sigma}(h(x))=\{1\}$ for $h(x)\in \mathcal{M}_{c_\sigma}$. Thus we deduce from $P_{c_\sigma}(h(x))=P_{c_\tau}((h(x))\pi)$ for $h(x)\in \mathcal{M}_{c_\sigma}$ that $s(\sigma)=id$ if and only if $s(\tau)=id$. Hence (i) and (ii) hold.

Now, it follows from (i) and (ii) that (1) and (2) are obviously true for the case $r(\sigma)=id$ or $s(\sigma)=id$.

\textbf{From now on, we further assume} both $r(\sigma)\ne id$ and $s(\sigma)\neq id$, and therefore $r(\tau)\ne id$ and $s(\tau)\neq id$ by (i) and (ii).

By the descriptions of $\mathcal{M}_{c_{r(\sigma)}}$ and $\mathcal{M}_{c_{r(\tau)}}$, the restriction of $\pi$ to $\mathcal{M}_{c_{r(\sigma)}}$ is mapped surjectively to $\mathcal{M}_{c_{r(\tau)}}$.
For $v(x)\in \mathcal{M}_{c_{r(\sigma)}}$, there holds $P_{c_{r(\sigma)}}(v(x))=\{1\}=P_{c_{r(\tau)}}((v(x))\pi)$. Thus $c_{r(\sigma)}$ and $c_{r(\tau)}$ are $M$-equivalent, and therefore $D$-equivalent by Corollary \ref{derp}(2).

In the sequel, we show that $c_{s(\sigma)}$ and $c_{s(\tau)}$ are $D$-equivalent, or equivalently, $M$-equivalent.

Actually, due to $s(\sigma)\neq id$, $\mathcal{M}_{c_{s(\sigma)}}$ consists of all reducible polynomials in $\mathcal{M}_{c_{\sigma}}$. By (ii), $\mathcal{M}_{c_{s(\tau)}}$ consists of all reducible polynomials in $\mathcal{M}_{c_{\tau}}$.
By the first condition of Definition \ref{newequrel}(1), the map $\pi$ sends irreducible polynomials to irreducible polynomials.
Thus the restriction of $\pi$ to $\mathcal{M}_{c_{s(\sigma)}}$ gives rise to a bijection between $\mathcal{M}_{c_{s(\sigma)}}$ and $\mathcal{M}_{c_{s(\tau)}}$.

By $(\gamma)$, there are positive integers $a, b$ such that $(x-1)^{p^a}\in \mathcal{M}_{c_{s(\sigma)}}$ and $(x-1)^{p^b}\in \mathcal{M}_{c_{s(\tau)}}$. We consider the two possible cases.

Case 1. $((x-1)^{p^a})\pi=(x-1)^{p^b}$. Then, by $(\delta)$, we have $$P_{c_{s(\sigma)}}((x-1)^{p^a}) = P_{c_\sigma}((x-1)^{p^a}) = P_{c_\tau}((x-1)^{p^b}) = P_{c_{s(\tau)}}((x-1)^{p^b}),$$ $P_{c_{s(\sigma)}}(h(x)) = P_{c_\sigma}(h(x))\setminus\{1\}$ for all $h(x)\in \mathcal{M}_{c_{s(\sigma)}}\setminus\{(x-1)^{p^a}\}$ and $P_{c_{s(\tau)}}(g(x))$ = $P_{c_\tau}(g(x))\setminus\{1\}$ for all $g(x)\in \mathcal{M}_{c_{s(\tau)}}\setminus\{(x-1)^{p^b}\}$. Thus, for $h(x)\in \mathcal{M}_{c_{s(\sigma)}}\setminus\{(x-1)^{p^a}\}$, the equality holds
$$P_{c_{s(\sigma)}}(h(x))=P_{c_\sigma}(h(x))\setminus\{1\}=P_{c_\tau}(h(x)\pi)\setminus\{1\}=P_{c_{s(\tau)}}(h(x)\pi).$$ This implies that the restriction of $\pi$ to $\mathcal{M}_{c_{s(\sigma)}}$ gives rise to an $M$-equivalence between $c_{s(\sigma)}$ and $c_{s(\tau)}$.

Case 2. $((x-1)^{p^a})\pi\neq (x-1)^{p^b}$. By the definition of $\pi$, we have an algebra isomorphism $R[x]/((x-1)^{p^a})\simeq R[x]/(((x-1)^{p^a})\pi)$. This implies that  $((x-1)^{p^a})\pi=(x+u)^{p^a}$ for some $u\in R$. Similarly, we may suppose $((x-1)^{p^b})\pi^{-1}=(x+v)^{p^b}$ for some $v\in R$. Due to $((x-1)^{p^a})\pi\neq (x-1)^{p^b}$, we have $u\ne -1$ and $v\ne -1$. Now we define a map $$\pi': \mathcal{M}_{c_{s(\sigma)}}\lra \mathcal{M}_{c_{s(\tau)}},$$
$$(x-1)^{p^a}\mapsto (x-1)^{p^b}, \; (x+v)^{p^b}\mapsto (x+u)^{p^a}, \; h(x)\mapsto (h(x))\pi \; $$for $h(x)\in \mathcal{M}_{c_{s(\sigma)}}\setminus \{(x-1)^{p^a},(x+v)^{p^a}\}.$ Then it follows from the bijection of $\pi$ that $\pi'$ is also a bijection.

We show that $\pi'$ defines an $M$-equivalence between $c_{s(\sigma)}$ and $c_{s(\tau)}$. By definition, it only remains to show that the corresponding power index sets are equal. In fact, by $(\delta)$, for $h(x)\in \mathcal{M}_{c_{s(\sigma)}}\setminus \{(x-1)^{p^a},(x+v)^{p^a}\}$, we have $P_{c_{s(\sigma)}}(h(x))=P_{c_{s(\tau)}}(h(x)\pi')$. So, to complete the proof, we have to show  $$P_{c_{s(\sigma)}}((x-1)^{p^a})=P_{c_{s(\tau)}}((x-1)^{p^b}) \mbox{ and } P_{c_{s(\sigma)}}((x+v)^{p^b})=P_{c_{s(\tau)}}((x+v)^{p^a}).$$On the one hand, $P_{c_{\sigma}}((x+v)^{p^b})\subseteq P_{c_{\sigma}}((x-1)^{p^a})$ by $(\alpha)$. Similarly, $P_{c_{\tau}}((x+u)^{p^a})\subseteq P_{c_{\tau}}((x-1)^{p^b})$.  On the other hand, by the definition of $M$-equivalences, we have $P_{c_{\sigma}}((x+v)^{p^b})=P_{c_{\tau}}((x-1)^{p^b})$ and $P_{c_{\sigma}}((x-1)^{p^a})=P_{c_{\tau}}((x+u)^{p^a})$. Thus $a=b$ and $$P_{c_{\sigma}}((x+v)^{p^b})=P_{c_{\tau}}((x-1)^{p^b})= P_{c_{\sigma}}((x-1)^{p^a})=P_{c_{\tau}}((x+u)^{p^a}).$$ Therefore it follows from $(\delta)$ that $P_{c_{s(\sigma)}}((x+v)^{p^b})=P_{c_{\sigma}}((x+v)^{p^b})\setminus\{1\}= P_{c_{\tau}}((x+u)^{p^a})\setminus\{1\}= P_{c_{s(\tau)}}((x+u)^{p^a})$ and $P_{c_{s(\sigma)}}((x-1)^{p^a})=P_{c_{\sigma}}((x-1)^{p^a})=P_{c_{\tau}}((x-1)^{p^b})=P_{c_{s(\tau)}}((x-1)^{p^b})$. Thus $c_{s(\sigma)}$ and $c_{s(\tau)}$ are $M$-equivalent. $\square$

\medskip
Generally, the converse of Proposition \ref{regular-singular} may be false, see Example \ref{ex4.5} in the next section.

\section{Examples and further questions\label{examples}}
In this subsection, we provide examples to illustrate results mentioned in the previous sections, and propose a few open questions for further considerations.

\begin{Bsp}\label{ex4.5} {\rm Let $R$ be an algebraically closed field of characteristic $5$. We take $\sigma\in \Sigma_{19}$ with the cycle type $(15,4)$, and $\tau\in \Sigma_{20}$ with the cycle type $(15,3,2)$. In this case, $r(\sigma)$ is a permutation of the cycle type $(4,1^{15})$ and $s(\sigma)$ is a permutation of cycle type $(15, 1^4)$, while $r(\tau)$ has the cycle type $(3,2, 1^{15})$ and $s(\tau)$ has the cycle type $(15,1^5)$. Clearly,  $S_{19}(c_{s(\sigma)},R)$ and $S_{20}(c_{s(\tau)},R)$ are derived equivalent by Corollary \ref{one-more}. Since
$\mathcal{M}_{c_{r(\sigma)}}= \{x-1, x+1, x-\eta, x+\eta\}$ and $\mathcal{M}_{c_{r(\tau)}}= \{x-1, x+1, x+\epsilon, x-\epsilon^2\}$, where $\eta$ and $\epsilon$ are $4$-th and $3$-th primitive roots of unity, respectively, it follows from Theorem \ref{main1} that
$S_{19}(c_{r(\sigma)},R)$ and $S_{20}(c_{r(\tau)},R)$ are derived equivalent.

By Lemma \ref{per}, $\mathcal{M}_{c_\sigma}=\{(x-1)^5,(x-\epsilon)^5,(x-\epsilon^2)^5,x+1,x-\eta,x+\eta\}$ and $\mathcal{M}_{c_\tau}=\{(x-1)^5,(x-\epsilon)^5,(x-\epsilon^2)^5,x+1\}$. Clearly, $|\mathcal{M}_{c_\sigma}|=6\neq 4=|\mathcal{M}_{c_\tau}|$. Hence there are no bijections between $\mathcal{M}_{c_\sigma}$ and $\mathcal{M}_{c_\tau}$, and therefore $S_{19}(c_{\sigma},R)$ and $S_{20}(c_{\tau},R)$ cannot be derived equivalent by Theorem \ref{main1}.

This shows that, in general, derived equivalences for both $p$-regular parts and $p$-singular parts of permutations do not have to guarantee a derived equivalence for the permutations themselves.
}\end{Bsp}

The following example shows that the existence of a Morita equivalence between centralizer matrix algebras depends on the ground field.

\begin{Bsp} {\rm Let $\sigma:=(1~2~3~4~5)(6~7~8\cdots~17~18), \tau:=(1~2~3~4~5~6~7)(8~9\cdots ~17~18)\in \Sigma_{18}$. The minimal polynomials of $c_\sigma$ and $c_\tau$ over $\mathbb{Q}$ are $(x-1)(x^4+x^3+x^2+x+1)(x^{12}+x^{11}+\cdots +x+1)$ and $(x-1)(x^{10}+x^9+\cdots+x+1)(x^6+x^5+\cdots +x+1)$, respectively. Moreover, $\mathcal{M}_{c_{\sigma}}=\{x-1, x^4+x^3+x^2+x+1, x^{12}+x^{11}+\cdots +x+1\}$ and $\mathcal{M}_{c_{\tau}}=\{x-1, x^{10}+x^9+\cdots+x+1, x^6+x^5+\cdots +x+1\}$. Clearly, there is no bijection between  $\mathcal{M}_{c_{\sigma}}$ and $\mathcal{M}_{c_{\tau}}$ such that all quotient algebras in Definition \ref{newequrel} (1) are isomorphic. Hence, by Theorem $\ref{main1}$, $S_{18}(c_\sigma, \mathbb{Q})$ and $S_{18}(c_\tau, \mathbb{Q})$ are not Morita equivalent, while $S_{18}(c_\sigma, \mathbb{C})$ and $S_{18}(c_\tau, \mathbb{C})$ are Morita equivalent (see also \cite[Theorem 1.2(2)]{xz2}). By Corollary \ref{derp}(2), this example also shows that derived equivalences of centralizer matrix algebras over $R$ depend upon the ground field $R$.
}\end{Bsp}

We point out that even in the class of centralizer matrix algebras, derived equivalences do not have to preserve representation-finiteness, while almost $\nu$-stable derived equivalences always preserve representation-finiteness for arbitrary algebras.

\begin{Bsp}{\rm
Let $R$ be an algebraically closed field, $c:=J_5(0)\oplus J_4(0)\oplus J_1(0)\in M_{10}(R)$ and $d:=J_5(0)\oplus J_2(0)\oplus J_1(0)\in M_8(R)$. Then $S_{10}(c,R)$ and $S_8(d,R)$ are derived equivalent by Theorem \ref{main1}, while $S_{10}(c,R)$ is representation-finite, but $S_8(d,R)$ is not by Lemma \ref{rep-f}.
}\end{Bsp}

Having described derived equivalences of centralizer matrix algebras, we propose the following questions for further study.

{\bf Question 1}. Let $R$ be a field. Under which conditions on permutations $\sigma\in S_n$ and $\tau\in S_m$ does the converse of Proposition \ref{regular-singular} hold true?

\smallskip
{\bf Question 2}. Let $R$ be a field and $c\in M_n(R)$. Is there any canonical form of the matrix $c$ under the equivalence relations in Definition \ref{newequrel}?

\smallskip
Related to generalization of Theorem \ref{main1} (see also Remark \ref{rmk3.3}), we mention the following.

{\bf Question 3}. Can one extend Theorem \ref{main1} to the case that $R$ is a principal ideal domain?

\medskip
\textbf{Acknowledgements.} The research work was supported partially by the National Natural Science Foundation of China (Grants 12031014 and 12226314).

{\footnotesize
\smallskip
Xiaogang Li,
School of Mathematical Sciences, Capital Normal University, 100048
Beijing, P. R. China; and

Shenzhen International Center for Mathematics, Southern University of Science and Technology, 518055
Shenzhen, Guangdong, P. R. China

{\tt Email: 2200501002@cnu.edu.cn}

\smallskip
Changchang Xi,
School of Mathematical Sciences, Capital Normal University, 100048
Beijing, P. R. China; and

School of Mathematics and Statistics, Shaanxi Normal University, 710119 Xi'an, P. R. China

{\tt Email: xicc@cnu.edu.cn}
}

\end{document}